\documentclass{article} 
\usepackage[a4paper,top=3cm,bottom=2cm,left=3cm,right=3cm,marginparwidth=1.75cm]{geometry}

\usepackage{amsmath, amssymb, amsthm}
\usepackage{mathtools}
\usepackage{bm}
\usepackage{booktabs}
\usepackage{multicol}
\usepackage{enumitem}
\usepackage[colorlinks=true, allcolors=blue]{hyperref}
\usepackage{blkarray}
\usepackage{xcolor}

\theoremstyle{definition}
\newtheorem{definition}{Definition}[section]
\newtheorem{remark}[definition]{Remark}
\newtheorem{example}[definition]{Example}

\theoremstyle{plain}
\newtheorem{theorem}[definition]{Theorem}
\newtheorem{lemma}[definition]{Lemma}

\newtheorem{corollary}[definition]{Corollary}

\newtheorem{prop}[definition]{Proposition}

\renewenvironment{proof}[1][\noindent Proof]{{\par\pushQED{\qed}\itshape #1\@. }}{\popQED}

\DeclareMathOperator{\I}{\mathrm{Im}}
\DeclareMathOperator{\PG}{PG}
\DeclareMathOperator{\PGL}{PGL}

\renewcommand{\L}{\mathcal{L}}
\renewcommand{\P}{\mathcal{P}}
\newcommand{\B}{\mathcal{B}}

\newcommand{\cQ}{\mathcal{Q}}
\newcommand{\cH}{\mathcal{H}}
\newcommand{\cL}{\mathcal{L}}

\newcommand{\N}{\mathbb{N}}
\newcommand{\F}{\mathbb{F}}

\newcommand{\bbR}{\mathbb{R}}

\newcommand{\comments}[1]{}
\newcommand{\gs}[3]{\genfrac{[}{]}{0pt}{}{#1}{#2}_{#3}}

\newcommand{\qbin}[2]{\genfrac{[}{]}{0pt}{}{#1}{#2}_{q}}
\allowdisplaybreaks

\title{Regular ovoids and Cameron-Liebler sets of generators in polar spaces}
\author{
	Maarten De Boeck\footnote{Department of Mathematical Sciences, University of Memphis, \href{mailto:mdeboeck@memphis.edu}{mdeboeck@memphis.edu} and Department of Mathematics: Algebra and Geometry, Ghent University},
	Jozefien D'haeseleer\footnote{Department of Mathematics: Analysis, Logic and Discrete Mathematics, Ghent University, Belgium, \href{mailto:jozefien.dhaeseleer@ugent.be}{jozefien.dhaeseleer@ugent.be}},
	Morgan Rodgers\footnote{Department of Mathematics, RPTU Kaiserslautern-Landau, Kaiserslautern, Germany,		\href{mailto:morgan.rodgers@rptu.de}{morgan.rodgers@rptu.de}}
}
\date{}

\begin{document}
	
	\maketitle

\begin{abstract}
	Cameron-Liebler sets of generators in polar spaces were introduced a few years ago as natural generalisations of the Cameron-Liebler sets of subspaces in projective spaces.
	In this article we present the first two constructions of non-trivial Cameron-Liebler sets of generators in polar spaces.
	Also regular $m$-ovoids of $k$-spaces are introduced as a generalization of $m$-ovoids of polar spaces.
	They are used in one of the aforementioned constructions of Cameron-Liebler sets.
\end{abstract}
	

	
	\section{Introduction}\label{sec:intro}
	
	In 1982, Cameron and Liebler introduced particular line sets in $\PG(n,q)$ when investigating the orbits of the action of $\PGL(n+1,q)$ on the points and lines of the projective space $\PG(n,q)$ (\cite{cameronliebler}).
	These line classes were later called Cameron-Liebler line sets in their honor.
	Many characterisation and classification results about Cameron-Liebler sets were obtained (see~\cite{bd, cossidentepavese, debeulemannaert, TF:14, gmp, gmet, gmol, met, met2, rod} amongst others).
	The many equivalent ways to describe a Cameron-Liebler set, both algebraically and combinatorially, sparked the interest of many researchers, and allowed for generalisations.
	Cameron-Liebler sets of $k$-spaces in $\PG(n,q)$ were studied in~\cite{bdd,bddc,debeulemannaertstorme,met3,rodstovan},
	and Cameron-Liebler sets of subsets of finite sets were discussed in~\cite{deboeckstormesvob, filmussets, meysets}.
	In~\cite{CLpol} Cameron-Liebler sets of generators in finite classical polar spaces were introduced, and those were further investigated in~\cite{CLpolequiv}.
	{More recently, Cameron-Liebler sets were also studied in the context of affine-classical spaces~\cite{CLaffineclassical}.}
	From a graph theory point of view Cameron-Liebler sets can be considered as intriguing sets in the sense of~\cite{debruynsuzuki}.
	The finite set, projective geometry and polar space context correspond to the Johnson, Grassmann and dual polar graphs.
	In many contexts Cameron-Liebler sets correspond to \emph{Boolean degree one functions} of the schemes corresponding to these graphs.
	\par In each of the contexts, equivalent characterisations of Cameron-Liebler sets exists.
	The central question is to classify the Cameron-Liebler sets.
	Typically, results are known for small parameters. Only for the set case a complete classification is known.
	The best results for subspaces of projective geometries are obtained in the case of line classes in $\PG(3,q)$.
	For projective geometries this is also the only case were non-trivial examples are known. For polar spaces so far only trivial examples were described.
	\par In this article we will present two constructions of non-trivial Cameron-Liebler sets of generators in polar spaces.
	Section~\ref{sec:prelim} deals with the necessary background of polar spaces, algebraic combinatorics and Cameron-Liebler sets.
	In Section~\ref{sec:exampleKleinquadric} we present a family of non-trivial Cameron-Liebler sets of generators in the hyperbolic quadrics $\mathcal{Q}^{+}(5,q)$, $q$ odd, and in Section~\ref{sec:newexample} we introduce and discuss a new construction of Cameron-Liebler sets for polar spaces with $e\leq 1$.
	The latter is based on a generalization of ovoids to higher-dimensional subspaces of polar spaces, which is discussed in Section~\ref{sec:ovoids}.
	In Section~\ref{sec:alcom} we investigate the eigenvalues of
	the distance-one relations on subspaces of polar spaces, with the lengthy calculations deferred to Appendix~\ref{ap:calculations}.
	These eigenvalue calculations are essential for the discussion of our generalization of ovoids.

	\section{Preliminaries}\label{sec:prelim}
	
	In this section we will first define polar spaces, then discuss distance-regular graphs and finally provide some background on Cameron-Liebler sets for polar spaces
	
	\subsection{Polar spaces}
	Polar spaces of rank $d\geq2$ are a kind of incidence geometries, whose axiomatic definition goes back to Veldkamp and Tits~\cite{tits,Veldkamp}.
	A polar space of rank $d$ has subspaces of dimension $1,\dots,d$. Polar spaces of rank 2 are called \emph{generalised quadrangles}.
	For reasons that we will not discuss in detail here, they are introduced separately.
	Points are considered subspaces of dimension 1.
	
	\begin{definition}
		A \emph{generalised quadrangle of order $(s,t)$}, with $s,t\geq1$, is a point-line geometry satisfying the following axioms.
		\begin{itemize}
			\item Two distinct points are incident with at most one line.
			\item Any line is incident with $s+1$ points and any point is incident with $t+1$ lines.
			\item For any point $P$ and line $\ell$, with $P$ not incident with $\ell$, there is a unique tuple $(P',\ell')$ such that $P$ is incident with $\ell'$, $P'$ is incident with $\ell$ and $P'$ is incident with $\ell'$.
		\end{itemize}
	\end{definition}
	
	\begin{definition}\label{def:combinatorialpolarspaces}
		A \emph{polar space of rank $d$}, $d\geq3$, is an incidence geometry $(\Pi,\Omega)$ with $\Pi$ a set whose elements are called points and $\Omega$ a set of subsets of $\Pi$ satisfying the following axioms.
		\begin{enumerate}
			\item Any element $\omega\in\Omega$ together with the elements of $\Omega$ that are contained in $\omega$, is a projective geometry of (algebraic) dimension at most $d$.
			\item The intersection of two elements of $\Omega$ is an element of $\Omega$ (the set $\Omega$ is closed under intersections).
			\item For a point $P\in\Pi$ and an element $\omega\in\Omega$ of dimension $d$ such that $P$ is not contained in $\omega$ there is a unique element $\omega'\in\Omega$ of dimension $d$ containing $P$ such that $\omega\cap\omega'$ is a hyperplane of $\omega$.
			The element $\omega$ is the union of all 2-dimensional elements of $\Omega$ that contain $P$ and are contained in $\omega$.
			\item There exist two elements $\Omega$ both of dimension $d$ whose intersection is empty.
		\end{enumerate}
	\end{definition}
	
	We now introduce the finite classical polar spaces.
	
	\begin{definition}
		A \emph{finite classical polar space} is an incidence geometry consisting of the totally isotropic subspaces of a non-degenerate quadratic or non-degenerate reflexive sesquilinear form on a vector space $\F^{n}_{q}$.
	\end{definition}

	A bilinear form for which all vectors are isotropic is called \emph{symplectic}; if $f(v,w)={f(w,v)}$ for all $v,w\in V$, then the bilinear form is called symmetric.
	A sesquilinear form on $V$ is called \emph{Hermitian} if the corresponding field automorphism $\theta$ is an involution and $f(v,w)={f(w,v)}^{\theta}$ for all $v,w\in V$.
	At various places in this article we will consider the finite classical  polar spaces as substructures of a projective space, in which they can naturally be embedded.
	We will however always use the \emph{vector space dimension}, which is one higher than the \emph{projective dimension}.
	The subspaces of dimension $1$ (vector lines), $2$ (vector planes) and $3$ (vector solids) are called \emph{points}, \emph{lines} and \emph{planes}, respectively. We now list the finite classical polar spaces of rank $d$.

	\begin{itemize}
		\item The hyperbolic quadric $\mathcal{Q}^{+}(2d-1,q)$ embedded in $\PG(2d-1,q)$. It arises from a hyperbolic quadratic form on $V(2d,q)$. Its standard equation is $X_{0}X_{1}+\dots+X_{2d-2}X_{2d-1}=0$.
		\item The parabolic quadric $\mathcal{Q}(2d,q)$ embedded in $\PG(2d,q)$. It arises from a parabolic quadratic form on $V(2d+1,q)$. Its standard equation is $X^{2}_{0}+X_{1}X_{2}+\dots+X_{2d-1}X_{2d}=0$.
		\item The elliptic quadric $\mathcal{Q}^{-}(2d+1,q)$ embedded in $\PG(2d+1,q)$. It arises from an elliptic quadratic form on $V(2d+2,q)$. Its standard equation is $g(X_{0},X_{1})+\dots+X_{2d-2}X_{2d-1}+X_{2d}X_{2d+1}=0$ with $g$ a homogeneous irreducible quadratic polynomial over $\F_{q}$.
		\item The Hermitian polar space $\mathcal{H}(2d-1,q)$ embedded in $\PG(2d-1,q)$ (where $q$ is a square). It arises from a Hermitian form on $V(2d,q)$, constructed using the field automorphism $x\mapsto x^{\sqrt{q}}$. Its standard equation is $X^{\sqrt{q}+1}_{0}+X^{\sqrt{q}+1}_{1}+\dots+X^{\sqrt{q}+1}_{2d-1}=0$.
		\item The Hermitian polar space $\mathcal{H}(2d,q)$ embedded in $\PG(2d,q)$ (where $q$ is square). It arises from a Hermitian form on $V(2d+1,q)$, constructed using the field automorphism $x\mapsto x^{\sqrt{q}}$. Its standard equation is $X^{\sqrt{q}+1}_{0}+X^{\sqrt{q}+1}_{1}+\dots+X^{\sqrt{q}+1}_{2d}=0$.
		\item The symplectic polar space $\mathcal{W}(2d-1,q)$ embedded in $\PG(2d-1,q)$. It arises from a symplectic form on $V(2d,q)$. For this symplectic form we can choose an appropriate basis $\{e_{1},\dots,e_{d},e'_{1},\dots,e'_{d}\}$ of $V(2d,q)$ such that $f(e_{i},e_{j})=f(e'_{i},e'_{j})=0$ and $f(e_{i},e'_{j})=\delta_{i,j}$, with $1\leq i,j\leq d$.
	\end{itemize}
	
	\begin{remark}\label{rem:tits}
		The fundamental result by Tits (see~\cite{tits}) states that all finite polar spaces of rank at least $3$ are finite classical polar spaces.
	\end{remark}
	
	In this paper, all polar spaces we will discuss are finite classical polar spaces, so we will refer to them as polar spaces.
	
	\begin{definition}
		The subspaces of maximal dimension (being $d$) of a polar space of rank $d$ are called \emph{generators}. We define the \emph{parameter} $e$ of a polar space $\mathcal{P}$ over $\F_{q}$ as $\log_{q}(x-1)$ with $x$ the number of generators through a $(d-1)$-space of $\mathcal{P}$.
	\end{definition}
	
	The parameter of a polar space only depends on the type of the polar space and not on its rank. Table~\ref{parameter} gives an overview.
	
	\begin{table}[ht]
		\centering
		\begin{tabular}{lc}\toprule
			polar space               & $e$ \\ \midrule
			$\mathcal{Q}^{+}(2d-1,q)$ & 0   \\
			$\mathcal{H}(2d-1,q)$     & 1/2 \\
			$\mathcal{W}(2d-1,q)$     & 1   \\
			$\mathcal{Q}(2d,q)$       & 1   \\
			$\mathcal{H}(2d,q)$       & 3/2 \\
			$\mathcal{Q}^{-}(2d+1,q)$ & 2   \\ \bottomrule
		\end{tabular}
		\caption{The parameter of the finite classical polar spaces.}\label{parameter}
	\end{table}
	
	An important concept, associated to polar spaces, are polarities.
	\begin{definition}
		A \emph{polarity} on $\PG(n,q)$ is an inclusion reversing involution $\perp$ acting on the subspaces of $\PG(n,q)$. In other words, $\perp^2$ is the identity, and any two subspaces $\pi$ and $\sigma$ satisfy $\pi \subseteq \sigma\Leftrightarrow\rho^\perp \subseteq \pi^\sigma$.
	\end{definition}
	
	Consider a non-degenerate sesquilinear form $f$ on the vector space $V=V(n+1,q)$, or the bilinear form $f$, based on a non-degenerate quadratic form $Q$ on the vector space $V=V(n+1,q)$, with $f(v,w)=Q(v+w)-Q(v)-Q(w)$. For a subspace $W$ of $V$, we can define its orthogonal complement with respect to $f$:
	\begin{align*}
		W^\perp=\{v\in V \ | \ \forall w\in W: f(v,w)=0\}.
	\end{align*}
	The  map $\perp$ that maps the subspace $W$ onto the  subspace $W^\perp$, is a polarity, and every polarity arises in this way.
	To every polar space a polarity is associated (but not the other way around).
	The image of a subspace $\pi$ on the polar space $\mathcal{P}$ under the corresponding polarity is its \emph{tangent space} $T_{\pi}(\mathcal{P})$, which is the subspace spanned by the lines meeting $\pi$ in a point and either fully contained in $\mathcal{P}$ or meeting $\mathcal{P}$ in precisely one point.
	\medskip
	
	We now introduce the \emph{Gaussian binomial coefficient}, which is very useful to describe counting results for vector spaces and polar spaces.
	
	\begin{definition}
		The \emph{Gaussian binomial coefficient} for integers $n,k$ with $n\geq k\geq0$ and prime power $q\geq2$ is given by
		\[
		\qbin{n}{k}=\prod^{k}_{i=1}\frac{q^{n-k+i}-1}{q^{i}-1}=\frac{(q^{n}-1)\cdots(q^{n-k+1}-1)}{(q^{k}-1)\cdots(q-1)}\:.
		\]
		For integers $n,k$ with $k>n\geq0$ or $n\geq0>k$ we set $\qbin{n}{k}=0$.
	\end{definition}
	
	The Gaussian binomial coefficient $\qbin{n}{k}$ equals the number of $k$-spaces in the vector space $\F^{n}_{q}$ (equivalently, the projective space $\PG(n-1,q)$).
	
	\begin{theorem}[{\cite[Theorem 9.4.1]{BCN} and~\cite[Remark 4.1.2]{phdvanhove}}]\label{thm:subspacecount}
		Given a classical polar space of rank $d$ and parameter $e$ defined over $\F_{q}$, and integers $a,b$ such that $0\leq b<a\leq d$, the number of $b$-spaces through a fixed $a$-space is given by
		\[
		\qbin{d-a}{b-a} \prod_{i=1}^{b-a} \left(q^{d-a+e-i} +1 \right)\:.
		\]
		The total number of $b$-spaces corresponds to the case $a=0$.
	\end{theorem}
	
	We also have the following useful lemma for Gaussian binomial coefficients, which is the $q$-analogue of the binomial theorem.
	
	\begin{lemma}\label{lem:gaussianbinomialidentity}
		For an integer $n\geq0$, prime power $q$ and $t$ any function independent of $n$, we have
		\[
		\sum_{k=0}^{n}\qbin{n}{k}q^{\binom{k}{2}}t^{k}=\prod_{i=0}^{n-1}(1+q^{i}t)\:.
		\]
		If $n=0$, the empty product in the right hand side is considered as 1.
	\end{lemma}

	\subsection{The algebraic combinatorics of polar spaces}
	
	In this section we will investigate the algebraic combinatorics underlying polar spaces, needed for the introduction of Cameron-Liebler sets in the next subsection;
	it will also allow us to discuss generalizations of ovoids to higher-dimensional subspaces in Section~\ref{sec:ovoids}.
	We refer to~\cite[Chapter 2]{BCN} and~\cite{del1} for a more extensive discussion of this topic.
	\par Throughout this paper, we will denote the all-ones $(n\times n)$-matrix by $J_{n}$, and the $(n\times n)$ identity matrix by $I_{n}$.
	The all-one vector of length $m$ will be denoted by $\bm{j_{m}}$. In all three cases the subscript will be omitted if the dimension is clear from the context. The first important definition is that of an association scheme.
	Note that in our work, we only consider association schemes in which all relations are symmetric.
	
	\begin{definition}
		A \emph{$d$-class association scheme} $\mathcal{R}$ on a finite non-empty set $\Omega$ is a set $\{R_{0},\dots,R_{d}\}$ of $d+1$ symmetric relations on $\Omega$ such that the following properties hold.
		\begin{enumerate}
			\item $R_{0}$ is the identity relation.
			\item $\mathcal{R}$ is a partition of $\Omega^{2}$.
			\item For $x,z\in\Omega$ with $(x,z)\in R_{k}$, the number of $y\in\Omega$ such that $(x,y)\in R_{i}$ and $(y,z)\in R_{j}$, equals a constant $p^{k}_{ij}$ (and is thus independent of $x$ and $z$).
		\end{enumerate}
	\end{definition}

	\par Each symmetric relation gives rise to its \emph{adjacency matrix}.
	This is the matrix whose rows and columns are indexed by the elements of the set on which the relation is defined and such that the entry on position $(x,y)$ equals 1 if the pair $(x,y)$ is in the relation and equals 0 otherwise.
	If $\mathcal{R}=\{R_{0},\dots,R_{d}\}$ is an association scheme on $\Omega$, with corresponding adjacency matrices $A_{i}$, then we have that $A_{0}=I$, all the $A_{i}$'s are symmetric, that $\sum^{d}_{i=0}A_{i}=J$ and that $A_{i}A_{j}=\sum^{d}_{k=0}p^{k}_{ij}A_{k}$.
	
	\begin{remark}\label{rem:alsographs}
		We point out that every antireflexive, symmetric relation on a set corresponds to a simple graph. In particular, the $d$ relations different from the identity, of a $d$-class association scheme on a set $\Omega$ correspond to a set of $d$ graphs that partition the complete graph on $\Omega$.
	\end{remark}
	
	\par We now give details on the association schemes arising from the classical polar spaces.
	This section follows the details from Chapter 4 of Vanhove's thesis~\cite{phdvanhove}.
	Recall that all dimensions given are vector-space dimensions.
	Throughout this subsection $\mathcal{P}$ is a classical polar space of rank $d$ and parameter $e$, defined over $\F_{q}$, and $G$ is its full automorphism group.
	We denote the collection of totally isotropic $a$-spaces of $\mathcal{P}$ by $\Omega_{a}$.
	Unlike for projective spaces (Grassmann scheme), the orbit under $G$ of an ordered pair $(\pi_{a}, \pi_{b}) \in \Omega_{a} \times \Omega_{b}$ depends not only on the intersection of the two subspaces, but also on the Witt index of $\langle \pi_{a}, \pi_{b} \rangle$.
	In fact we have the following.
	
	\begin{theorem}[{\cite[Proposition 4.9]{Stanton2}}]\label{th:relations}
		The orbits of $G$ on $\Omega_{a} \times \Omega_{b}$ are given by
		\[
		R_{a,b}^{s,k} = \left\lbrace \left(\pi_{a}, \pi_{b}\right) \mid \dim{\left( \pi_{a} \cap \pi_{b}\right)} = s, \dim\left( \langle \pi_{a}, \pi_{b} \cap \pi_{a}^{\perp} \rangle \right)=k \right\rbrace,
		\]
		with $0 \leq s \leq \min{(a,b)}$, $\max{(a,b)} \leq k \leq \min{(d, a+b-s)}$.
	\end{theorem}
	
	The relation $R_{a,b}^{\min{(a,b)}, \max{(a,b)}}$ just defines incidence between the subspaces (or equality if $a=b$) and we denote this relation by $R_{a,b}$.
	Note that for any fixed $a$, the relations $R_{a,a}^{s,k}$ define an association scheme on $\Omega_{a}$; in this case the relations with $s=a-1$ are called \emph{distance-one relations}.

	We write $C_{a,b}^{s,k}$ for the $(0,1)$-matrix associated with the relation $R_{a,b}^{s,k}$ (the rows are labelled by the $b$-spaces and the columns by the $a$-spaces); in particular we write $C_{a,b}$ for the $(0,1)$-matrix associated with the relation $R_{a,b}$.
	We remark that each $C_{a,b}^{s,k}$ can be thought of as giving a homomorphism $\bbR^{\Omega_{a}} \to \bbR^{\Omega_{b}}$.
	For convenience of notation in some upcoming results, we will consider $C_{a,b}^{s,k} = 0$ if there is no corresponding (non-empty) relation $R_{a,b}^{s,k}$ on $\Omega_{a} \times \Omega_{b}$.

	We now give an important result from~\cite{phdvanhove}, which will be used often in the following sections.
	Note that $\bbR^{\Omega_{a}}$ is the $\bbR$-vector space whose positions are indexed by the element of $\Omega_{a}$.
	\begin{theorem}\label{thm:PolarSpaces}
		Let $\mathcal{P}$ be a 
		polar space of rank $d$.
		\begin{enumerate}[label=(\roman*)]
			\item Under the action of $G$, the module $\bbR^{\Omega_{a}}$ has a unique decomposition into orthogonal irreducible subspaces given by
			\[
			\bbR^{\Omega_{a}} = \bigoplus\limits_{\substack{0 \leq r \leq a\\0 \leq i \leq \min{(r,d-a)}}} V_{r,i}^{a}
			\]
			where
			\[
			V_{r,i}^{a} = \I{(C_{r,a})} \cap \I{(C_{r-1,a})}^{\perp} \cap \I{(C_{d-i,a})} \cap \I{(C_{d-i+1,a})}^{\perp}\:.
			\]
			In particular for each $a$, $V_{0,0}^{a} = \langle \bm{j} \rangle$, where $\bm{j}\in\bbR^{\Omega_{a}}$ is the all-ones vector.
			\item The submodules $V_{r,i}^{a} \subseteq \bbR^{\Omega{a}}$ and $V_{r,i}^{b} \subseteq \bbR^{\Omega_{b}}$ are isomorphic when both are defined.
			\item\label{thmpart:submoduleimage} The restriction of the incidence map $C_{a,b} :  \bbR^{\Omega_{a}} \to \bbR^{\Omega_{b}}$ to the submodule $V_{r,i}^{a}$ is trivial when $V_{r,i}^{b}$ is not defined, and is a bijection between the isomorphic submodules otherwise.
			\item The restriction of every map $C_{a,b}^{s,k}$ to $V_{r,i}^{a}$ is a scalar multiple of the restriction of $C_{a,b}$.
		\end{enumerate}
	\end{theorem}
	The fact that these $V_{r,i}^{a}$ are irreducible subspaces under the action of $G$ gives us the fact that each eigenspace of each $C_{a,a}^{s,k}$ is a direct sum of the $V_{r,i}^{a}$.
	
	\begin{definition}\label{def:combinatorialdesign}
		We say that a set $M \subseteq \Omega_{a}$ is a \emph{combinatorial design with respect to the $b$-spaces} if, for any $\pi_{b} \in \Omega_{b}$, the size of the set $\{\pi_{a} \in M \mid (\pi_{a}, \pi_{b}) \in R_{a,b}^{s,k} \}$ is a constant depending only on $s$ and $k$.
	\end{definition}
	
	For the special case of a combinatorial design with respect to the generators of $\mathcal{P}$, this condition simplifies to: for every $\pi_{b} \in \Omega_{b}$, the number of elements $\pi_{a} \in M$ meeting $\pi_{b}$ in an $s$-space is a constant depending only on $s$. We have the following characterizations due to Ito.
	
	\begin{theorem}\label{thm:combinatorialdesign}
		Let $M \subseteq \Omega_{a}$. Then $M$ is a combinatorial design with respect to the $b$-spaces of $\mathcal{P}$ if and only if its characteristic vector $\bm{\chi}$ is orthogonal to every irreducible submodule $V^{a}_{r,i}$ in $\bbR^{\Omega_{a}}$ which has an isomorphic copy in $\bbR^{\Omega_{b}}$ different from $V_{0,0}^{a}$.
	\end{theorem}

	\subsection{Cameron-Liebler sets for polar spaces}
	
	We consider the set of generators $\Omega_{d}$ of a 
	polar space
	of rank $d$, and the corresponding vector space $\bbR^{\Omega_{d}}$. The \emph{incidence vector} of a set $S\subseteq\Omega_{d}$ is the vector in $\bbR^{\Omega_{d}}$ having a 1 at every position that corresponds to an element of $S$ and a 0 elsewhere.
	On $\Omega_{d}$ we have the relations $R^{i,d}_{d,d}=R_{i}$ for $i \in \{ 0,\dots,d \}$, corresponding to the matrices $C_{i}=C^{i,d}_{d,d}$, which can be interpreted as linear maps on $\bbR^{\Omega_{d}}$.
	Recall that  $\mathclap{R} \ =\{R_{0},R_{1},\dots,R_{d}\}$ is a $d$-class association scheme. In view of Remark~\ref{rem:alsographs} the graph corresponding to $R_1$ is called the \emph{dual polar graph} of $\mathcal{P}$.
	This is a so-called \emph{distance-regular graph}. Each relation in $\mathcal{R}$ can also be seen as a distance relation in the dual polar graph of $\mathcal{P}$.
	\par We also know that
	\[
	\bbR^{\Omega_{d}} = \bigoplus_{r=0}^{d}V_{r}\:,
	\]
	where $V_{r}=V_{r,0}^{d}=\I{(C_{r,d})} \cap \I{(C_{r-1,d})}^{\perp}$ and thus $V_0=\langle\bm{j}\rangle$. As we mentioned before the eigenspaces of the matrices $C_{0},\dots,C_{d}$ are direct sums of the subspaces $V_{r}$, $r=0,\dots,d$.
	Consequently, each subspace $V_r$ corresponds to an eigenvalue of $C_i$, but it might happen that several $V_r$'s correspond to the same eigenvalue of $C_i$ for some $i$. We now introduce (degree one) Cameron-Liebler sets following the definitions given in~\cite{CLpolequiv, CLpol}.
	
	\begin{definition}\label{def:CLset}
		Let $\mathcal{L}\subseteq\Omega_{d}$ be a set of generators of a polar space of rank $d$ with parameter $e$ defined over $\F_{q}$, with incidence vector $\chi_{\mathcal{L}}$.
		Then $\mathcal{L}$ is a \emph{degree one Cameron-Liebler space} if $\chi_{\mathcal{L}}\in V_0\oplus V_1$. It is called a \emph{Cameron-Liebler set} if $\chi_{\mathcal{L}}\in V_0\oplus V$ where $V$ is the eigenspace of $C_{d}$ to which $V_1$ belongs.
		The \emph{parameter} of a (degree one) Cameron-Liebler set is
		\[
		x=\frac{|\mathcal{L}|}{\prod_{i=0}^{d-2}(q^{e+i}+1)}\:.
		\]
	\end{definition}
	
	For most polar spaces the definition of a Cameron-Liebler set and a degree one Cameron-Liebler set coincide; the exceptions are the hyperbolic quadrics of even rank, the parabolic quadrics of odd rank en the symplectic polar spaces of odd rank, i.e.~the polar spaces $\mathcal{Q}^{+}(4n-1,q)$, $\mathcal{Q}(4n+2,q)$ and $\mathcal{W}^{+}(4n+1,q)$ for all $q$ and $n$.
	See~\cite[Section 3]{CLpol} for a detailed discussion. Clearly, each degree one Cameron-Liebler set is a Cameron-Liebler set.
	
	\begin{remark}\label{rem:equivalent}
		In~\cite{CLpolequiv,CLpol} several equivalent definitions for Cameron-Liebler sets are presented. From the description above it is clear that $S$ is a degree one Cameron-Liebler set if and only if its incidence vector $\chi_{S}$ is contained in $\I(C_{1})$, or in other words, if and only if there exists a vector $v\in\bbR^{\Omega_1}$ such that $C_{1,d}v=\chi_{S}$.
		We also have the following result.
	\end{remark}
	
	\begin{theorem}[{\cite[Theorem 3.1]{CLpolequiv}}]\label{stellinggg}
		Let $\mathcal{P}$ be a 
		polar space of rank $d$ with parameter $e$ defined over $\F_{q}$, let $\mathcal{L}$ be a set of generators of $\mathcal{P}$ and let $i$ be an integer with $1\leq i\leq d$.
		If $\mathcal{L}$ is a degree one Cameron-Liebler set of generators in $\mathcal{P}$, with parameter $x$,  then the number of elements of $\mathcal{L}$ meeting a generator $\pi$ in a $(d-i)$-space equals
		\begin{align}\label{formulelang}
			\begin{dcases}
				\left( (x-1)\begin{bmatrix} d-1 \\ i-1 \end{bmatrix} + q^{i+e-1}\begin{bmatrix} d-1 \\ i \end{bmatrix}\right)q^{\binom{i-1}{2}+ (i-1)e}
				& \text{if $\pi \in \mathcal{L}$}     \\
				\hfil x \begin{bmatrix} d-1 \\ i-1\end{bmatrix} q^{\binom{i-1}{2}+(i-1)e}
				& \text{if $\pi \notin \mathcal{L}$.}
			\end{dcases}
		\end{align}
		Moreover, if this property holds for a polar space $\mathcal{P}$ and an integer $i$ such that
		\begin{itemize}
			\item $i$ is odd if $\mathcal{P}=Q^+(2d-1,q)$,
			\item $i\neq d$ if $d$ is odd, and $\mathcal{P}=Q(2d,q)$ or $\mathcal{P}=W(2d-1,q)$
			\item $i$ is arbitrary if $\mathcal{P}$ is another polar space,
		\end{itemize}
		then $\mathcal{L}$ is a degree one Cameron-Liebler set with parameter $x$.
	\end{theorem}
	
	We present now two examples of (degree one) Cameron-Liebler sets.
	
	\begin{example}[{\cite[Example 4.2]{CLpol}}]\label{ex:pointpencil}
		A \emph{point-pencil} with \emph{vertex} $P$ in a polar space $\mathcal{P}$ is the set of all generators in $\mathcal{P}$ containing the point $P$. A point-pencil is a degree-one Cameron-Liebler set (and thus also a Cameron-Liebler set) since its incidence vector is clearly in $\I(C_{1})$. It has parameter 1.
		\par A \emph{partial ovoid} is a set of points on $\mathcal{P}$ which are pairwise not collinear. In other words, for any two points the line through them is not isotropic with respect to the underlying quadratic or sesquilinear form.
		Hence for any two points of a partial ovoid, the point-pencils with these vertices are disjoint. So, If $\mathcal{O}$ is a partial ovoid, the union of the point-pencils with the points of $\mathcal{O}$ as vertex, is a (degree one) Cameron-Liebler set with parameter $|\mathcal{O}|$.
	\end{example}
	
	\begin{example}[{\cite[Example 4.4]{CLpol}}]\label{ex:embedded}
		Let $\mathcal{P}$ be a 
		polar space of rank $d$ with parameter $e$ over $\F_{q}$ and let $\mathcal{P}'$ be a 
		polar space of the same rank with parameter $e-1$ over $\F_{q}$ that is embedded in $\mathcal{P}$, $e\geq1$; all examples of such polar spaces $\mathcal{P}$ and $\mathcal{P}'$ are listed below.
		Let $\mathcal{L}$ be the set of generators contained in $\mathcal{P}'$, note that these are generators of both $\mathcal{P}$ and $\mathcal{P}'$.
		Then $\mathcal{L}$ is a Cameron-Liebler set with parameter $q^{e-1}+1$.
		\par We present all examples that arise from this construction. Up to a projection (symplectic polar space) they are
		\begin{itemize}
			\item There are parabolic quadrics $\mathcal{Q}(2d,q)$ embedded in the elliptic quadric $\mathcal{P}=\mathcal{Q}^{-}(2d+1,q)$. Each of them gives rise to a (degree one) Cameron-Liebler set of $\mathcal{P}$ with parameter $q+1$.
			\item There are hyperbolic quadrics $\mathcal{Q}^{+}(2d-1,q)$ embedded in the parabolic quadric $\mathcal{P}=\mathcal{Q}(2d,q)$. Each of them gives rise to a (degree one) Cameron-Liebler set of $\mathcal{P}$ with parameter $2$.
			Recall that the symplectic variety $\mathcal{W}(2d-1,q)$ is isomorphic to $\mathcal{P}$ if $q$ is even.
			\item There are Hermitian polar spaces $\mathcal{H}(2d-1,q)$ embedded in the Hermitian polar space $\mathcal{P}=\mathcal{H}(2d,q)$, $q$ a square.
			Each of them gives rise to a (degree one) Cameron-Liebler set of $\mathcal{P}$ with parameter $\sqrt{q}+1$.
		\end{itemize}
	\end{example}
	
	We also have the following.
	
	\begin{lemma}[{\cite[Lemma 3.3]{CLpolequiv} and~\cite[Lemma 4.1]{CLpol}}]\label{basicoperations}
		Let $\mathcal{P}$ be a 
		polar space. Let $\mathcal{L}$ and $\mathcal{L}'$ be (degree one) Cameron-Liebler sets of $\mathcal{P}$ with parameters $x$ and $x'$ respectively.
		\begin{itemize}
			\item[(i)] $0\leq x\leq q^{e+d-1}+1$.
			\item[(ii)] The complement of $\mathcal{L}$ in $\Omega_{d}$ is a (degree one) Cameron-Liebler set of $\mathcal{P}$ with parameter $q^{e+d-1}+1-x$.
			\item[(iii)] If $\mathcal{L}\cap\mathcal{L}'=\emptyset$, then $\mathcal{L}\cup\mathcal{L}'$ is a (degree one) Cameron-Liebler set of $\mathcal{P}$ with parameter $x+x'$.
			\item[(iv)] If $\mathcal{L}'\subset\mathcal{L}$, then $\mathcal{L}\setminus\mathcal{L}'$ is a (degree one) Cameron-Liebler set of $\mathcal{P}$ with parameter $x-x'$.
		\end{itemize}
	\end{lemma}
	
	The (degree one) Cameron-Liebler sets from Examples~\ref{ex:pointpencil} and~\ref{ex:embedded}, their complements and disjoint unions, are called \emph{trivial}; all others are called non-trivial.
	All degree one Cameron-Liebler sets known so far are trivial examples. There are examples of Cameron-Liebler sets that are not degree one Cameron-Liebler sets, see~\cite[Table 4]{CLpolequiv}.
	\par Finally we mention a few classification results. Note that from Definition~\ref{def:CLset} it follows that the parameter of a (degree one) Cameron-Liebler set is in $\mathbb{Q}$, but actually we have the following.
	
	\begin{theorem}[{\cite[Lemma 4.8]{CLpol}} and {\cite[Lemma 5.3]{CLpolequiv}}]
		Let $\mathcal{P}$ be a 
		polar space. If $\mathcal{L}$ is a (degree one) Cameron-Liebler set of $\mathcal{P}$ with parameter $x$, then $x\in\N$.
	\end{theorem}
	
	The most important classification result so far states that small (degree one) Cameron-Liebler sets are trivial.
	
	\begin{theorem}[{\cite[Theorem 5.5]{CLpolequiv}}]\label{th:classificationsmalldegreeoneCL}
		Let $\mathcal{P}$ be a 
		polar space of rank $d$ and parameter $e$, and let $\mathcal{L}$ be a degree one Cameron-Liebler set of $\mathcal{P}$ with parameter $x$.
		If $x\leq q^{e-1}+1$, then $\mathcal{L}$ is the union of $x$ point-pencils whose vertices are pairwise non-collinear or $x=q^{e-1}+1$ and  $\mathcal{L}$ is the set of generators in an embedded polar space of rank $d$ and with parameter $e-1$.
	\end{theorem}
	
	\begin{theorem}[{\cite[Theorem 6.7]{CLpol}}]\label{th:classificationsmallCL}
		Let $\mathcal{P}$ be a 
		polar space of rank $d$ and parameter $e$, which is not a hyperbolic quadric of even rank, a parabolic quadric of odd rank or a symplectic polar space of odd rank.
		Let $\mathcal{L}$ be a Cameron-Liebler set of $\mathcal{P}$ with parameter $x$. If $x\leq q^{e-1}+1$, then $\mathcal{L}$ is the union of $x$ point-pencils whose vertices are pairwise non-collinear or $x=q^{e-1}+1$ and $\mathcal{L}$ is the set of generators in an embedded polar space of rank $d$ and with parameter $e-1$.
	\end{theorem}
	
	From Theorem~\ref{th:classificationsmalldegreeoneCL} it follows that all degree one Cameron-Liebler sets with parameter 1 are point-pencils. This does not generalise to Cameron-Liebler sets.
	In~\cite[Theorem 6.4]{CLpol} Cameron-Liebler sets with parameter 1 that are not trivial are described, showing that the exclusion of several types of polar spaces in Theorem~\ref{th:classificationsmallCL} is necessary.
	\par For the symplectic polar space $\mathcal{W}(5,q)$ and the parabolic quadric $\mathcal{Q}(6,q)$, it is proven in~\cite[Theorem 5.9]{CLpolequiv} that a CL set with parameter $x$, with $2\leq x\leq \sqrt[3]{2q^2}-\frac{\sqrt[3]{4q}}{3}+\frac{1}{6}$ is a union of embedded polar spaces $Q^+(5,q)$ and point-pencils.

\section{The example in the Klein quadric \texorpdfstring{$\mathcal{Q}^+(5,q)$}{Q+(5,q)}}\label{sec:exampleKleinquadric}

In this section we provide a family of new examples of non-trivial degree one Cameron-Liebler sets of generators in $\mathcal{Q}^+(5,q)$, $q$ odd.
As mentioned in the previous section, we know that for these polar spaces, the parameter of the CL set is an integer between $0$ and $q^2+1$, and CL sets with parameter $1$ are point-pencils.
Since $\cQ^+(5,q)$ contains an ovoid of size $q^2+1$, namely an embedded elliptic quadric $\cQ^-(3,q)$, we know that for every parameter $x\in \{0, \dots, q^2+1\}$, there exists a trivial CL example, namely the union of $x$ disjoint point-pencils.
This family of examples was first described in the PhD thesis of the second author~\cite{PhDDhaeseleer}, and are the first known nontrivial CL sets in $Q^+(5,q)$.
To give the construction, we use the Klein correspondence $\kappa$ between the lines of $\PG(3,q)$  and the points of the Klein quadric $\mathcal{K}=\mathcal{Q}^+(5,q)$~\cite{Klein}.

Recall that the generators of a hyperbolic quadric $\mathcal{Q}^+(2n+1,q)$ can be divided in two classes, commonly called the Latin and Greek generators. By the Klein correspondence, the points of a Latin plane in $\mathcal{K}$ correspond to the set of lines through a fixed point in $\PG(3,q)$ (called a \emph{point-pencil}, the fixed point is called the \emph{vertex}), and the points of a Greek plane in $\mathcal{K}$ correspond to the set of lines in a fixed plane in $\PG(3,q)$.

Consider the hyperbolic quadric $Q\cong \mathcal{Q}^+(3,q)$ in $\PG(3,q)$, defined by the equation $X_0 X_1+X_2 X_3=0$. Via $\kappa$ the lines of $Q$ correspond to the set of points of two conics $C\cup C'$ in $\mathcal{K}$. Moreover, the planes $\alpha = \langle C \rangle$ and $\alpha'=\langle C' \rangle$, are each others image under the polarity corresponding to $\mathcal{K}$.

Every point $P\in \PG(3,q)$ gives rise to a Latin plane $\pi^P_l$ and a Greek plane $\pi^P_g$ in $\mathcal{K}$: $\pi^P_l$ and $\pi^P_g$ are the images under $\kappa$ of the set of lines through $P$ and the set of lines in the plane $P^\perp$, respectively. Here $\perp$ is the polarity related to the quadric $Q$ in $\PG(3,q)$.

\begin{definition}
	A point $P(x_0,x_1,x_2,x_3) \in \PG(3,q)$ is a \emph{square point} if $x_0 x_1+x_2 x_3$ is a non-zero square in $\mathbb{F}_q$. A point $P(x_0,x_1,x_2,x_3) \in \PG(3,q)$ is a \emph{non-square point} if $x_0 x_1+x_2 x_3$ is a non-square in $\mathbb{F}_q$.
\end{definition}

Now we can partition the set of planes in $\mathcal{K}$ into the following sets.
\begin{multicols}{2}
	\begin{itemize}
		\item $\mathcal{S}_l= \{\pi^P_l | P $ is a square point$ \}$
		\item $\mathcal{NS}_l= \{\pi^P_l | P $ is a non-square point$ \}$
		\item $\mathcal{O}_l= \{\pi^P_l | P\in Q \}$
	\end{itemize}
	\begin{itemize}
		\item $\mathcal{S}_g= \{\pi^P_g | P $ is a square point$ \}$
		\item $\mathcal{NS}_g= \{\pi^P_g |  P $ is a non-square point$ \}$
		\item $\mathcal{O}_g= \{\pi^P_g | P\in Q \}$
	\end{itemize}
\end{multicols}
It is known that a $2$-secant to $Q$ in $\PG(3,q)$, $q$ odd, contains $\frac{q-1}{2}$ square points and $\frac{q-1}{2}$ non-square points. A line disjoint from $Q$ in $\PG(3,q)$ contains $\frac{q+1}{2}$ square points and $\frac{q+1}{2}$ non-square points. For a tangent line $\ell$ to $Q$, there are two possibilities; $\ell$ contains $q$ square points, or $\ell$ contains $q$ non-square points, see~\cite[Table 15.5(c)]{Hirschfeld3}. In the first case, $\ell$ is a \emph{square tangent line}. In the latter case, $\ell$ is a \emph{non-square tangent line}.

We partition the set of points in $\mathcal{K}$ into the following sets.
\begin{itemize}
	\item The set $\mathcal{X}_{1S}$ of points in $\mathcal{K}$ corresponding to the square tangent lines to $Q$.
	\item The set $\mathcal{X}_{1NS}$ of points in $\mathcal{K}$ corresponding to the non-square tangent lines to $Q$.
	\item The set $\mathcal{X}_{2}$ of points in $\mathcal{K}$ corresponding to the $2$-secants to $Q$.
	\item The set $\mathcal{X}_{0}$ of points in $\mathcal{K}$ corresponding to the lines disjoint from $Q$.
	\item The set $\mathcal{X}_{\infty}=C\cup C'$ of points in $\mathcal{K}$ corresponding to the lines of $Q$.
\end{itemize}
We present two lemmas that will be useful in the remainder of the construction.

\begin{lemma}\label{lemmarechte}
	If $l$ is a square tangent line to $Q$ in $\PG(3,q)$, then $l^\perp$ is a square tangent line if $q\equiv 1 \mod 4$, and $l^\perp$ is a non-square tangent line if $q\equiv 3 \mod 4$. 	If $l$ is a non-square tangent line to $Q$ in $\PG(3,q)$, then $l^\perp$ is a non-square tangent line if $q\equiv 1 \mod 4$, and $l^\perp$ is a square tangent line if $q\equiv 3 \mod 4$.
\end{lemma}
\begin{proof}
	Consider a tangent line $l$ to $Q$ in $\PG(3,q)$. Since the orthogonal group $PGO_+(4,q)$ of $\mathcal{Q}^+(3,q)$ acts transitively on the points of $Q=\mathcal{Q}^+(3,q)$ (see~\cite[Theorem 22.6.4]{thashirschfeld}), we may suppose that $l$ contains the point $(1,0,0,0)$ of $Q$,
	and so $l=\langle (1,0,0,0),(0,0,1,t)\rangle$, for a fixed $t\in \mathbb{F}_q\setminus \{0\}$. Note that $l$ is a square tangent line if and only if $t$ is a square in $\mathbb{F}_q$.
	We find that $T_{(1,0,0,0)}(Q)$ is the plane defined by $x_1=0$, while $T_{(0,0,1,t)}(Q)$ is the plane defined by $tX_2+X_3=0$. The intersection of these two planes is $l^\perp = \langle (1,0,0,0), (0,0, 1, -t) \rangle$. The lemma follows since $l^\perp$ is a square line if and only if $-t$ is a square in $\mathbb{F}_q$, and $-1$ is a square $\mathbb{F}_q$ if and only if $q\equiv 1 \mod 4$.
\end{proof}

The following lemma follows directly from~\cite[Theorem $22.7.2$]{thashirschfeld}.

\begin{lemma}\label{thasje}
	If $l$ is a bisecant to $Q$ in $\PG(3,q)$, then $l^\perp$ is also a bisecant to $Q$. Furthermore, if $l$ is a line skew to $Q$ in $\PG(3,q)$, then $l^\perp$ is also skew to $Q$.
\end{lemma}

We now continue with the construction of the non-trivial Cameron-Liebler set. We will use the theory of tactical decompositions (see~\cite[p. 7]{tactdec}) to prove that the constructed set is a Cameron-Liebler set.

\begin{definition}
	Let $(\P, \B, I)$ be an incidence geometry with $\P$ a set of points and $\B$ a set of blocks. Let $\{ P_1, P_2, \dots , P_s  \}, P_i\neq \emptyset$ be a partition of $\P$, and let $\{ B_1, B_2, \dots , B_r  \}, B_i\neq \emptyset$ be a partition of $\B$.
	\begin{itemize}
		\item If there exists an $(s\times r)-$matrix $X$ with $|\{p\in P_i| \ p \ I \ b\}|=X_{ij}, \forall b\in B_j$, then the decomposition is called \emph{block-tactical}\index{block-tactical decomposition}.
		\item If there exists an $(s\times r)-$matrix $Y$ with $|\{b\in B_i| \ p \ I \ b\}|=Y_{ij}, \forall p\in P_j$, then the decomposition is called \emph{point-tactical}\index{point-tactical decomposition}.
	\end{itemize}
	The decomposition is called \emph{tactical}\index{tactical decomposition} if it is both block- and point-tactical.
\end{definition}

\begin{lemma}\label{lemmatactical}
	Let $(\P, \B, I)$ be an incidence geometry with $\P$ a set of points, $\B$ a set of blocks and $A$ the point-block incidence matrix. Let $\{ P_1, P_2, \dots , P_s  \}, P_i\neq \emptyset$ be a partition of $\P$, and let $\{ B_1, B_2, \dots , B_r  \}, B_i\neq \emptyset$ be a partition of $\B$.
	\begin{itemize}
		\item If the partition is block-tactical with corresponding matrix $X$, then  \begin{align*}
			A^T \chi_{\P_i}=\sum_{l=1}^r X_{il}\chi_{\B_l}, \forall i\in \{1,\dots, s\}.
		\end{align*}
		\item If the partition is point-tactical with corresponding matrix $Y$, then  \begin{align*}
			A \chi_{\B_i}=\sum_{l=1}^s Y_{lj}\chi_{\P_l}, \forall i\in \{1,\dots, r\}.
		\end{align*}
	\end{itemize}
\end{lemma}

In the following proposition, we prove that the point and plane partitions
\begin{gather*}
	\{\mathcal{X}_{1S},\mathcal{X}_{1NS}, \mathcal{X}_{2}, \mathcal{X}_0, \mathcal{X}_{\infty}\} \\
	\{\mathcal{S}_l,\mathcal{S}_g,\mathcal{NS}_l,\mathcal{NS}_g,\mathcal{O}_l,\mathcal{O}_g   \}
\end{gather*}
give a point-tactical decomposition.

\begin{prop}\label{propvb}
	The partition of the points
	$\{\mathcal{X}_{1S},\mathcal{X}_{1NS}, \mathcal{X}_{2}, \mathcal{X}_{0}, \mathcal{X}_{\infty}\}$
	and the partition of the planes
	$\{\mathcal{S}_l,\mathcal{S}_g,\mathcal{NS}_l,\mathcal{NS}_g,\mathcal{O}_l,\mathcal{O}_g   \}$ of $\mathcal{K}$
	give a point-tactical decomposition with matrix $B_1$ if $q \equiv 1 \mod 4$ and the matrix $B_3$ if $q \equiv 3 \mod 4$.
	\begin{align*}
		B_1 = 	\begin{blockarray}{ccccccc}
			{\color{gray}\mathcal{S}_l} & {\color{gray} \mathcal{S}_g} & {\color{gray}\mathcal{NS}_l} &{\color{gray} \mathcal{NS}_g} & {\color{gray}\mathcal{O}_l} &{\color{gray}\mathcal{O}_g} \\
			\begin{block}{(cccccc)c}
				q&q&0&0&1&1 & {\color{gray}\mathcal{X}_{1S}} \\ 0&0&q&q&1&1& {\color{gray}\mathcal{X}_{1NS}}\\ \frac{q-1}{2}&\frac{q-1}{2}&\frac{q-1}{2}&\frac{q-1}{2}&2&2 &{\color{gray}\mathcal{X}_{2}} \\
				\frac{q+1}{2}&\frac{q+1}{2}&\frac{q+1}{2}&\frac{q+1}{2}&0&0 &{\color{gray}\mathcal{X}_{0}} \\
				0&0&0&0&q+1&q+1& {\color{gray}\mathcal{X}_{\infty}}\\
			\end{block}
		\end{blockarray}
	\end{align*}
	\[
	B_3 = 	\begin{blockarray}{ccccccc}
		{\color{gray}\mathcal{S}_l} & { \color{gray}\mathcal{S}_g} & {\color{gray} \mathcal{NS}_l} & { \color{gray}\mathcal{NS}_g} &  {\color{gray}\mathcal{O}_l} & {\color{gray}\mathcal{O}_g} \\
		\begin{block}{(cccccc)c}
			q&0&0&q&1&1 &  {\color{gray}\mathcal{X}_{1S}} \\ 0&q&q&0&1&1&  {\color{gray}\mathcal{X}_{1NS}}\\ \frac{q-1}{2}&\frac{q-1}{2}&\frac{q-1}{2}&\frac{q-1}{2}&2&2 & {\color{gray}\mathcal{X}_{2}} \\
			\frac{q+1}{2}&\frac{q+1}{2}&\frac{q+1}{2}&\frac{q+1}{2}&0&0 & {\color{gray}\mathcal{X}_{0} }\\
			0&0&0&0&q+1&q+1&  {\color{gray}\mathcal{X}_{\infty}}\\
		\end{block}
	\end{blockarray}
	\]
\end{prop}
\begin{proof}
	We find these matrices using the Klein correspondence, and so we will prove the lemma for the lines of $\PG(3,q)$ instead of the points of $\mathcal{K}$.
	This means that we will use point-pencils of lines and the lines in fixed planes of $\PG(3,q)$, instead of the planes in $\mathcal{K}$.
	
	Let $\ell$ be a line corresponding to the type of the row, so a square tangent line, non-square tangent line, bisecant, skew line or line contained in $Q$.
	The first and third element of the row correspond to the number of square and non-square point on $\ell$, respectively, since this is also the number of point-pencils with square and non-square vertex.
	The fifth and sixth element of the row correspond to the number of points of $Q$ on $\ell$. To compute the second and fourth element on the row, we recall that $\ell\subset R^\perp \iff R\in \ell^\perp$, with $R\in \PG(3,q)$.
	So, the second and fourth element equal the first and the third element on the row corresponding to $\ell^{\perp}$.
	Via Lemmas~\ref{lemmarechte} and~\ref{thasje} we know that a $\ell^{\perp}$ has the same type of $\ell$ unless $q\equiv3\mod{4}$ and $\ell$ is a square or non-square tangent line. Those two types are interchanged.
\end{proof}

\begin{theorem}\label{thmvb} Let $q$ be an odd prime power.
	\begin{itemize}
		\item The sets $\mathcal{S}_l\cup \mathcal{S}_g$, $\mathcal{NS}_l\cup \mathcal{NS}_g$ and $\mathcal{O}_l\cup \mathcal{O}_g$ are degree one Cameron-Liebler sets of planes in $\mathcal{K}$, with parameter $\frac{q(q-1)}{2}$, $\frac{q(q-1)}{2}$ and $q+1$ respectively, for $q\equiv 1 \mod 4$.
		\item The sets $\mathcal{S}_l\cup \mathcal{NS}_g$, $\mathcal{S}_g\cup \mathcal{NS}_l$ and $\mathcal{O}_l\cup \mathcal{O}_g$ are degree one Cameron-Liebler sets of planes in $\mathcal{K}$, with parameter $\frac{q(q-1)}{2}$, $\frac{q(q-1)}{2}$ and $q+1$ respectively, for $q\equiv 3 \mod 4$.
	\end{itemize}
\end{theorem}

\begin{proof}
	We prove this theorem for $q\equiv 3 \mod 4$. The proof for $q\equiv 1 \mod 4$ is analogous. Let $A$ be the point-plane incidence matrix of $\mathcal{K}$. We write $\chi_K$ for the incidence vector of a set of planes $K$, and likewise write $\chi_L$ for the incidence vector of a set $L$ of lines.
	From the previous proposition we find the following equations.
	\begin{align*}
		A^T\chi_{1S}     & =q\chi_{\mathcal{S}_l}+q\chi_{\mathcal{NS}_g}+\chi_{\mathcal{O}_l}+\chi_{\mathcal{O}_g}                                                             \\
		A^T\chi_{1NS}    & =q\chi_{\mathcal{S}_g}+q\chi_{\mathcal{NS}_l}+\chi_{\mathcal{O}_l}+\chi_{\mathcal{O}_g}                                                             \\
		A^T\chi_{2}      & =\frac{q-1}{2}( \chi_{\mathcal{S}_l}+\chi_{\mathcal{S}_g}+\chi_{\mathcal{NS}_l}+\chi_{\mathcal{NS}_g})+2(\chi_{\mathcal{O}_l}+\chi_{\mathcal{O}_g}) \\
		A^T\chi_{\infty} & =(q+1)(\chi_{\mathcal{O}_l}+\chi_{\mathcal{O}_g}).
	\end{align*}
	After some calculations, we find:
	\begin{align*}
		\chi_{\mathcal{S}_l}+\chi_{\mathcal{NS}_g} & =A^T\left( \frac{3q+1}{2q(q+1)}\chi_{1S}+\frac{q-1}{2q(q+1)}\chi_{1NS}-\frac{1}{q+1}\chi_{2}  \right) \\
		\chi_{\mathcal{S}_g}+\chi_{NS_l}           & =A^T\left( \frac{q-1}{2q(q+1)}\chi_{1S}+\frac{3q+1}{2q(q+1)}\chi_{1NS}-\frac{1}{q+1}\chi_{2}  \right) \\
		\chi_{\mathcal{O}_l}+\chi_{\mathcal{O}_g}  & = \frac{1}{q+1} A^T \chi_\infty.
	\end{align*}
	The sets $\mathcal{S}_l\cup \mathcal{NS}_g$, $\mathcal{S}_g\cup \mathcal{NS}_l$ and $\mathcal{O}_l\cup \mathcal{O}_g$ are contained in the image of $A^T$, and so they are degree one Cameron-Liebler sets of planes in $\mathcal{K}$, for $q\equiv 3 \mod 4$ by Remark~\ref{rem:equivalent}. The parameters of these Cameron-Liebler sets follow immediately from their size, see Definition~\ref{def:CLset}.
	
	Analogously, we find that the sets $\mathcal{S}_l\cup \mathcal{S}_g$, $\mathcal{NS}_l\cup \mathcal{NS}_g$ and $\mathcal{O}_l\cup \mathcal{O}_g$ are degree one Cameron-Liebler sets of planes in $\mathcal{K}$, for $q\equiv 1 \mod 4$.
\end{proof}

\begin{remark}
	Note that the Cameron-Liebler sets $\mathcal{O}_l\cup \mathcal{O}_g$ are the union of $q+1$ point-pencils, whose points are the elements of the conic $C$. Moreover, this set is also the set of point-pencils whose points are the elements of the conic $C'$. Hence, this example is a trivial Cameron-Liebler set. The other determined Cameron-Liebler sets in Theorem~\ref{thmvb} are new and non-trivial examples, in the sense that they are not a union of point-pencils. They even do not contain a single point-pencil.
\end{remark}

\begin{prop}\label{prop:Q+5qexamplenopointpencil}
	The sets $\mathcal{S}_l\cup \mathcal{S}_g$, and $\mathcal{NS}_l\cup \mathcal{NS}_g$, \ for $q\equiv 1 \mod 4$,
	and the  sets $\mathcal{S}_l\cup \mathcal{NS}_g$ and $\mathcal{S}_g\cup \mathcal{NS}_l$, \ for $q\equiv 3 \mod 4$
	{ do not contain a point-pencil. In particular, they }are not the union of point-pencils whose points are pairwise non-collinear.
\end{prop}
\begin{proof}
	We prove this proposition for the set $\L = \mathcal{S}_l\cup \mathcal{S}_g$, if $q\equiv 1 \mod 4$. The proofs for the other cases are analogous.
	{Suppose that $\L$ contains a point-pencil, and let $P$ be its vertex.}
	By investigating the sum of the first two columns of the matrix $B_1$ in Proposition~\ref{propvb}, we find that $P$ contains $2q$, $0$, $q-1$, $q+1$ or $0$ elements of $\L$ if it is contained in $\mathcal{X}_{1S}$, $\mathcal{X}_{1NS}$, $\mathcal{X}_{2}$, $\mathcal{X}_{0}$, or $\mathcal{X}_{\infty}$, respectively.
	In each case we find that $\L$ cannot contain all planes of $\mathcal{K}$ through $P$, a contradiction.
\end{proof}

\begin{remark}
	In Theorem~\ref{thmvb} two non-trivial Cameron-Liebler sets were described. We point out that there are isomorphisms of $\PG(3,q)$ that fix $Q$ and interchange the sets of square points and non-square points, e.g.~the isomorphism induced by the matrix $\nu I_{4}$ with $\nu$ a non-square in $\F_q$ and $I_{4}$ the identity matrix of order 4.
	These isomorphisms induce an isomorphism of $\mathcal{K}$ that interchanges $\mathcal{S}_l\cup \mathcal{S}_g$ and $\mathcal{NS}_l\cup \mathcal{NS}_g$, and also $\mathcal{S}_l\cup \mathcal{NS}_g$ and $\mathcal{S}_g\cup \mathcal{NS}_l$.
	So, Theorem~\ref{thmvb} describes for each $q$ a unique non-trivial Cameron-Liebler set of generators of $\mathcal{Q}^+(5,q)$, up to isomorphism.
\end{remark}

\begin{remark}
	The authors verified through computer search using Gurobi~\cite{Gur} that the Cameron-Liebler set described above is the only non-trivial Cameron-Liebler set of generators in $\cQ^+(5,3)$, up to isomorphism. There are no non-trivial Cameron-Liebler sets of generators in $\cQ^+(5,2)$.
\end{remark}


\section{Eigenvalues of polar spaces}\label{sec:alcom}

From Theorem~\ref{thm:PolarSpaces} and the subsequent remarks we know that a subspace $V^{a}_{r,i}$ is contained in an eigenspace of the map $C^{s,k}_{a,a}$,
and hence corresponds to an eigenvalue of $C^{s,k}_{a,a}$ (equivalently, its corresponding relation).
In this section, we are interested in calculating the values of some of these eigenvalues. We use the following notation.

\begin{definition}\label{def:mu}
	We denote the eigenvalue of $C^{s,k}_{n,n}$ on $V^{n}_{r,i}$ by $\mu_{d,n,s,k,r,i}$. Here we have $0\leq r\leq n$, $0\leq i\leq\min(r,d-n)$, $0\leq s\leq n$ and $n\leq k\leq\min(d,2n-s)$.
\end{definition}

A recursive approach to compute this eigenvalues exist; we will come back to it below.
But it is hard to find closed formulas (which allow to compare eigenvalues). To our knowledge such closed formulas were described for $n=d$~\cite[Theorem 5.4]{Stanton2}, $(s,k)=(n-1,\min(n+1,d))$~\cite[Section 4.3.2]{phdvanhove}, $(s,k)=(0,n)$~\cite[Section 4.3.3]{phdvanhove} and $(s,k,n)=(0,d-1,d-1)$~\cite[Lemma 13]{metsch},
where the parameters not mentioned are general (within the bounds allowed for by Definition~\ref{def:mu}).
In this section we will discuss these eigenvalues for some other specific choices of the parameters, namely $(s,k)=(n-1,n)$ and $(n,i)=(d-1,1)$.

We need the following counting results. The notation is based on~\cite{eisfeld}, and it follows from Theorem~\ref{th:relations} that it is well-defined.

\begin{definition}\label{def:alpha}
	For a 
	polar space of rank $d$ and with parameter $e$, defined over $\F_{q}$, the following values are well-defined.
	\begin{enumerate}
		\item For any $a$-space $\pi_{a}$, let $\alpha_{d,a,c,\ell}$ be the number of $c$-spaces $\pi_{c}$ with $(\pi_{a}, \pi_{c}) \in R_{a,c}^{0,\ell}$.
		\item For $\pi_{a}\in\Omega_{a}$ and $\pi_{b}\in\Omega_{b}$ with $(\pi_{a}, \pi_{b}) \in R_{a,b}^{s,k}$, let $\alpha_{d,(a,b),(s,k),c,(t,\ell)}$ be the number of $c$-spaces $\pi_{c}$ through $\pi_{b}$ with $(\pi_{a},\pi_{c}) \in R_{a,c}^{t,\ell}$.
	\end{enumerate}
\end{definition}

\begin{lemma}[{\cite[Lemma 9.4.2]{BCN} and~\cite[Theorem 3.6(d)]{eisfeld}}]\label{lem:alphavalues}
	For a 
	polar space of rank $d$ and parameter $e$, defined over $\F_{q}$ we have
	\begin{align*}
		\alpha_{d,a,c,\ell}                & = q^{\nu} \qbin{a}{\ell-c}\qbin{d-a}{\ell -a}      \prod_{i=0}^{\ell-a-1} \left( q^{d-l+e+i} +1 \right)\ \\
		\alpha_{d, (a,b),(s,k),c,(t,\ell)} & = \qbin{k-b}{t-s} \alpha_{d-b+s-t, k-b+s-t, c-b+s-t, \ell-b+s-t}
	\end{align*}
	where $\nu = a(\ell-a)+\frac{1}{2}(a+c-\ell)(a-c-3\ell+4d+2e-1)$.

\end{lemma}
Note that the definition of $\alpha_{d, (a,b),(s,k),c,(t,\ell)}$ depends on $a$, but the formula does not. We use this lemma for $b=0$, $s=0$, and $k=a$, to find the following corollary.

\begin{corollary}[{\cite[Lemma 9.4.2]{BCN}}]\label{countingBCN}
	For any $a$-space $\pi_a$ in a 
	polar space $\mathcal{P}$ of rank $d$ and parameter $e$, the number of $c$-spaces $\pi_c$ intersecting
	$\pi_a$ in an $t$-space and with Witt index $l$ for $\langle \pi_a, \pi_c\rangle$,
	i.e.\ with $(\pi_a, \pi_c)\in R_{a,c}^{t,l}$, is given by:
	\begin{align*}
		q^{\nu}
		\cdot \qbin{a}{t}\qbin{a-t}{l-c}\qbin{d-a}{l-a}\prod_{i=0}^{l-a-1}\left(q^{d-a-i+e-1}+1\right)
	\end{align*}
	where $\nu = (a-t)(l-a)+(a+c-t-l)(a-c-t-3l+4d+2e-1)/2$.
\end{corollary}

We also know from Theorem~\ref{thm:PolarSpaces}\ref{thmpart:submoduleimage} that the restriction of any map $C_{a,b}^{s,k}: \bbR^{\Omega_{a}} \to \bbR^{\Omega_{b}}$
to one of the $V_{r,i}^{a}$ is a scalar multiple of the restriction of $C_{a,b}$ to this space.
This allows us to introduce the following notation, inspired by~\cite[Definitions 4.3.1 and 4.3.2]{phdvanhove}. Compared to~\cite{phdvanhove}, we have introduced the parameter $d$ in the index of $\theta$, for sake of clarity.

\begin{definition}
	For a 
	polar space of rank $d$ we denote the unique scalar used in the restriction of $C_{a,b}^{s,k}$ to $V_{r,i}^{a}$ by $\theta_{d,(r,i),a,b,(s,k)}$. In other words for every $\bm{v} \in V_{r,i}^{a}$ we have
	\[
	C_{a,b}^{s,k}\bm{v} = \theta_{d,(r,i),a,b,(s,k)}C_{a,b}\bm{v}\:.
	\]
	Here we assume $0\leq r\leq d$, $0\leq r\leq \min(a,b)$, $0\leq i\leq\min(r,d-a,d-b)$, $0\leq s\leq\min(a,b)$ and $\max(a,b)\leq k\leq\min(d,a+b-s)$.
\end{definition}

The conditions on the parameters of $\theta_{d,(r,i),a,b,(s,k)}$ reflect that the relation $R^{s,k}_{a,b}$ and the subspace $V^{a}_{r,i}$ exist, and that $C_{a,b}$ is not the zero map on $V^{a}_{r,i}$
(which is by Theorem~\ref{thm:PolarSpaces}\ref{thmpart:submoduleimage} equivalent to the fact that $V^{b}_{r,i}$ exists). We define notation for two specific cases.

\begin{definition}\label{def:psichi}
	We define $\psi_{d,r,i,s,k}$ and $\chi_{d,i,r,n,t,\ell}$ as follows.
	\begin{enumerate}
		\item $\psi_{d,r,i,s,k} = \theta_{d,(r,i),r,d-i,(s,k)}$ with $0\leq s\leq r\leq d-i\leq k\leq\min(d,d-i+r-s)$ and $0\leq i\leq r$. 
		\item $\chi_{d,i,r,n,t,\ell} = \theta_{d,(r,i),r,n,(t,\ell)}$ with $0\leq t\leq r\leq n\leq \ell \leq \min(d,r+n-t)$ and $0\leq i\leq\min(r,d-n)$.
		We set $\chi_{d,i,r,n,t,\ell} =0$ if these conditions are not met.
	\end{enumerate}
\end{definition}

Although not made explicit in the notation, we remark that $\theta_{d,(r,i),a,b,(s,k)}$, $\psi_{d,r,i,s,k}$ and $\chi_{d,i,r,n,t,\ell}$ also depend on the parameter $e$ and the size of the underlying field $\F_{q}$ of the polar space.
The following result gives the value of $\chi_{d,i,r,n,t,\ell}$ for the case $n=\ell$ and $t=r-1$.

\begin{lemma}[{\cite[Lemma 4.3.9]{phdvanhove}}]\label{lem:psichivalues}

	For integers $d,r,i,n$ with $0 \leq i \leq r\leq n \leq d-i$ and $r \geq 1$ we have
	\[
	\chi_{d,i,r,n,r-1,n} = q^{d-i-n}\left(\qbin{i}{1}(q^{i-1+e}+1)- \qbin{r}{1}\right)\:.
	\]
\end{lemma}

\subsection{Eigenvalues related to the distance-one relations}

For $1 \leq a \leq d-1$ there are precisely two non-trivial distance-one relations on $a$-spaces, namely $R_{a,a}^{a-1,a}$ and $R_{a,a}^{a-1,a+1}$.
The latter relation gives the well-studied \emph{graph of Lie type} on $a$-spaces, and its eigenvalues are well known.

\begin{theorem}[{\cite[Theorem 4.3.10]{phdvanhove}}]\label{thm:eigenvaluesLie}
	Let $\mathcal{P}$ be a 
	polar space of rank $d$ and parameter $e$, defined over $\F_{q}$. Then the eigenvalue of $R_{a,a}^{a-1,a+1}$ on $V_{r,i}^{a}$, $1 \leq a \leq d-1$, is given by
	\begin{multline*}
		\qbin{a-r}{1}\qbin{d-a}{1} q \left(q^{d-a-1+e}+1\right) \\
		+ \qbin{a-r+1}{1} \left( \qbin{r}{1} \left(q^{d-a-i}-1\right) - \qbin{i}{1} q^{d-a-i}\left(q^{i-1+e}+1\right)\right)\:.
	\end{multline*}
\end{theorem}

In this subsection we investigate the eigenvalues of the former relation on the $V_{r,i}^{a}$, where $0 \leq r \leq a$, $0 \leq i \leq \min{(r,d-a)}$. To do this we follow the method given in~\cite[Section 4.3.2]{phdvanhove}.

We can follow the same method to show the following.

\begin{theorem}\label{thm:eigenvaluesNotLie}
	Let $\mathcal{P}$ be a 
	polar space of rank $d$ and parameter $e$, defined over $\F_{q}$. For each $a$ with  $1 \leq a \leq d-1$, the eigenvalue of $R_{a,a}^{a-1,a}$ on $V_{r,i}^{a}$ is given by
	\[
	\mu_{d,a,a-1,a,r,i}=\qbin{a-r}{1} q^{2d-2a+e} - \left( \qbin{r}{1} - \qbin{i}{1}(q^{i-1+e}+1) \right)q^{d-a-i}\:.
	\]
\end{theorem}
\begin{proof}
	Recall that $0\leq r\leq a$ and $0\leq i\leq\min(r,d-a)$. To consider the eigenvalue of $R_{a,a}^{a-1,a}$ on $V_{r,i}^{a}$, we look at the product $C_{a,a}^{a-1,a}C_{r,a}$.
	Consider a pair $(\pi_{r}, \pi_{a}) \in \Omega_{r}\times \Omega_{a}$.
	If there is an $a$-space $\pi_{a}^{\prime} \in \Omega_{a}$ such that $\pi_{r} \subseteq \pi_{a}^{\prime}$ -- in other words $(\pi_{r}, \pi_{a}^{\prime})\in R_{r,a}$ --
	and such that $(\pi_{a}^{\prime}, \pi_{a}) \in R_{a,a}^{a-1,a}$, then either $\pi_{r} \subseteq \pi_{a}$, or else $(\pi_{r}, \pi_{a}) \in R_{r,a}^{r-1,a}$.
	Using the notation from Definition~\ref{def:alpha} and then applying Lemma~\ref{lem:alphavalues} we have that
	\begin{align*}
		C_{a,a}^{a-1,a}C_{r,a} & = \alpha_{d,(a,r),(r,a),a,(a-1,a)}C_{r,a} + \alpha_{d,(a,r),(r-1,a),a,(a-1,a)}C_{r,a}^{r-1,a} \\
		& = \qbin{a-r}{1}\alpha_{d-a+1,1,1,1}C_{r,a} + \alpha_{d-a,0,0,0}C_{r,a}^{r-1,a}                \\
		& = \qbin{a-r}{1} q^{2d-2a+e} C_{r,a} + C_{r,a}^{r-1,a}.
	\end{align*}
	Now if $\bm{v} \in V_{r,i}^{a}$ then $\bm{v} = C_{r,a}\bm{w}$ for some $\bm{w} \in V_{r,i}^{r}$ by Theorem~\ref{thm:PolarSpaces}\ref{thmpart:submoduleimage}. Therefore
	\begin{align*}
		C_{a,a}^{a-1,a}\bm{v} & =C_{a,a}^{a-1,a}C_{r,a}\bm{w}                                                \\
		& =\qbin{a-r}{1}q^{2d-2a+1}\left(C_{r,a} \bm{w}\right) + C_{r,a}^{r-1,a}\bm{w} \\
		& =\qbin{a-r}{1}q^{2d-2a+1}\bm{v} + C_{r,a}^{r-1,a}\bm{w}\:.
	\end{align*}
	Now using the notation from Definition~\ref{def:psichi} and applying Lemma~\ref{lem:psichivalues} we see that
	\begin{align*}
		C_{r,a}^{r-1,a}\bm{w} & = \chi_{d,i,r,a,r-1,a} \left(C_{r,a}\bm{w}\right)     = \chi_{d,i,r,a,r-1,a} \bm{v}                                                                                                                                  & =q^{d-i-a} \left(\qbin{i}{1}(q^{i-1+e}+1) -\qbin{r}{1}\right) \bm{v}\:,
	\end{align*}
	giving the result.
\end{proof}

Given these eigenvalues we can prove the following result.

\begin{theorem}\label{th:uniqueeigenvalues}
	For a be a 
	polar space $\mathcal{P}$ the eigenvalues of the relation $R_{a,a}^{a-1,a}$ determine $V_{r,i}^{a}$ uniquely.
\end{theorem}
\begin{proof}
	The eigenvalue $\mu_{d,a,a-1,a,r,i}$ of $R_{a,a}^{a-1,a}$ on $V_{r,i}^{a}$ is given in Theorem~\ref{thm:eigenvaluesNotLie}, and we denote it by
	\[
	\lambda_{r,i}=\qbin{a-r}{1} q^{2d-2a+e} - \left( \qbin{r}{1} - \qbin{i}{1}(q^{i-1+e}+1) \right)q^{d-a-i}
	\]
	as $a$ and $d$ are fixed. We define
	\begin{align*}
		A_{r,i} & =\left((q-1)\lambda_{r,i}+q^{2d-2a+e}-q^{d-a}+q^{d-a-1+e}\right)q^{a-d} \\
		& =q^{d-r+e}-q^{r-i}+q^{i+e-1}\:.
	\end{align*}
	It is sufficient to prove that $R\to\bbR:(r,i)\mapsto A_{r,i}$ is injective, with $R=\{(r,i):0\leq r\leq a, 0\leq i \leq\min(r,d-a)\}$. Note that
	\[
	d_{r,i}=\max\{k:q^{k}\vert A_{r,i}\}=
	\begin{cases}
		i+e-1 & r>2i+e-1 \\
		d-r+e & r=2i+e-1 \\
		r-i   & r<2i+e-1
	\end{cases}
	\]
	and that the coefficient of $q^{d_{r,i}}$ in $A_{r,i}$ is 1 in the first two cases, and $-1$ in the final case. So, if $A_{r,i}=A_{r',i'}$ for some $(r,i),(r',i')\in R$, then one of the following cases occurs.
	\begin{itemize}
		\item We have $r>2i+e-1$, $r'>2i'+e-1$ and $i=i'$. From $q^{d-r+e}-q^{r-i}=q^{d-r'+e}-q^{r'-i}$ it then follows that also $r=r'$.
		\item We have $r>2i+e-1$, $r'=2i'+e-1$ and $i-1=d-r'$. We find that $q^{d-r+e}=q^{r-d+2i'+e-2}$, hence $i'=d-r+1$, a contradiction since $i'\leq d-a\leq d-r$.
		\item We have $r=2i+e-1$, $r'=2i'+e-1$ and $r=r'$. Immediately also $i=i'$.
		\item We have $r<2i+e-1$, $r'<2i'+e-1$ and $r-i=r'-i'$. We find that
		\begin{align*}
			q^{d-r+e}+q^{i+e-1}=q^{d-r'+e}+q^{i+r'-r+e-1}
		\end{align*}
		Without loss of generality we can assume that $r'>r$. We rewrite the previous equation as
		\begin{align*}
			q^{e-1}\left(q^{r'-r}-1\right)\left(q^{d-r'+1}-q^{i}\right)=0
		\end{align*}
		and we conclude that $r'=d-i+1$ and hence also $i'=d-r+1$. However, $i'\leq d-a\leq d-r'<d-r<d-r+1$, a contradiction.
	\end{itemize}
	So, in each case we find that $(r,i)=(r',i')$ or a contradiction. This concludes the proof.
\end{proof}

\subsection{Eigenvalues related to the \texorpdfstring{$(d-1)$}{(d-1)}-spaces}

In Sections~\ref{sec:ovoids} and~\ref{sec:newexample} we will have special attention for the $(d-1)$-spaces of a polar space. For that reason, we now focus on the eigenvalues in the association scheme corresponding to the $(d-1)$-spaces.
Most of the calculations are quite lengthy, so we defer them to Appendix~\ref{ap:calculations}. Here we will give an overview of the results.

The first main result is a general formula for the $\psi_{d,r,i,s,k}$.

\begin{lemma}
	For integers $d,r,i,s,k$ with $0\leq s\leq r\leq d-i\leq k\leq\min(d,d-i+r-s)$ and $0\leq i\leq r$ we have that
	\begin{align*}
		\psi_{d,r,i,s,k} & =
		\begin{multlined}[t]
			{(-1)}^{r-s}q^{\nu_{1}} \qbin{i}{d-k}  \prod_{m=1}^{\crampedclap{k-d+i}}(q^{i-m+e}+1)
			\sum_{\ell=0}^{\crampedclap{r-k+d-i-s}}{(-1)}^{\ell}\qbin{d-k}{\ell}\qbin{r-i}{k-d+s+\ell}q^{\nu_{2}}\:.
		\end{multlined}
	\end{align*}
	where
	\begin{align*}
		\nu_{1} & = \cramped{\binom{k-d+i}{2}+\binom{r-k+d-i-s}{2}+i(r-k+d-i-s)} \\
		\nu_{2} & = \ell(e-1-r+k-d+i+s+\ell)
	\end{align*}
\end{lemma}
\begin{proof}
	See Corollary~\ref{cor:psivalues}.
\end{proof}

The second result is a general formula for $\chi_{d,1,r,n,t,n+1}$.

\begin{lemma}
	For integers $d,r,n,t$ with $0\leq t\leq r\leq n \leq \min(d,r+n-t)-1$ and $r\geq 1$ we have that
	\begin{align*}
		\chi_{d,1,r,n,t,n+1} & = \begin{multlined}[t]{(-1)}^{r-t}q^{\binom{r-t}{2}+(r-t-1)(d-n-2)}\qbin{r-1}{t}\\ \cdot \left((q^{e}+1)\left(q^{d-t-2-n+r} \vphantom{\qbin{a}{b}} \right.\right.
			\left.\left. +(q^{r-t-1}-1)\qbin{d-n-1}{1}\right)
			-(q^{r}-1)\qbin{d-n-1}{1}\right).\end{multlined}
	\end{align*}
	
\end{lemma}
\begin{proof}
	See Lemma~\ref{lem:chi_n+1=l}.
\end{proof}

We use this result on the $\chi_{d,1,r,n,t,n+1}$ values to determine several eigenvalues of the scheme.

\begin{theorem}\label{th:eigenval}
	For integers $d,r,s,k$ with $1\leq r\leq d-1$, $0\leq s\leq d-1$ and $d-1\leq k\leq\min(d,2d-s-2)$ we have that the eigenvalue of $C^{s,k}_{d-1,d-1}$ on $V^{d-1}_{r,1}$ equals
	\begin{align*}
		& \mu_{d,d-1,s,k,r,1}=q^{e(d-r-s-1)}\sum^{s}_{t=0} \qbin{d-r-1}{s-t}{(-1)}^{r-t}q^{\binom{r-t}{2}+\binom{d-r-s+t+1}{2}+et-1}\:\mu'_{d,k,r,t} \\
		\text{with }\quad & \mu'_{d,k,r,t}=
		\begin{cases}
			\left(\qbin{r}{t}-(q^{e}+1)\qbin{r-1}{t}\right) & k=d-1 \\
			(q^{e}+1)\qbin{r-1}{t}                          & k=d
		\end{cases}
	\end{align*}
\end{theorem}
\begin{proof}
	See Theorem~\ref{th:cor:eigenval}.
\end{proof}

\section{\texorpdfstring{$(m,a)$}{(m,a)}-ovoids of 
	polar spaces}\label{sec:ovoids}

Among the classical substructures of polar spaces \emph{ovoids} are among the most studied. They were introduced in~\cite{thas} and we refer to~\cite{debeulestorme} for a recent overview.
An \emph{$m$-ovoid} is a set of points of the polar space such that each generator of the polar space contains precisely $m$ of these points.
In the particular case of the a 1-ovoid (shortly, ovoid), this means that no 2 points of the ovoid are collinear. We now generalise this concept.

\begin{definition}\label{def:maovoid}
	For a 
	polar space $\mathcal{P}$ of rank $d$, and $1\leq a\leq d-1$, we will call a set $M \subseteq \Omega_{a}$ an \emph{$(m,a)$-ovoid} of $\mathcal{P}$ if every generator of $\mathcal{P}$ contains precisely $m$ elements of $M$.
	\par An $(m,a)$-ovoid $M$ is called \emph{trivial} if either $M=\emptyset$ or $M=\Omega_{a}$.
\end{definition}

We can phrase the condition of an $(m,a)$-ovoid $M$ in terms of the matrices associated to the polar space: if $\bm{\chi}$ is the characteristic vector of $M$, then
\[
C_{a,d}\bm{\chi} = m \bm{j}\:.
\]
This implies that $\bm{\chi}$ is in the inverse image of $V_{0,0}^{d}$ under $C_{a,d}$ and so by Theorem~\ref{thm:PolarSpaces}\ref{thmpart:submoduleimage} we have that
\[
\bm{\chi} \in {\left( \bigoplus\limits_{1 \leq r \leq a} V_{r,0}^{a} \right)}^{\perp}=V_{0,0}^{a}\ \oplus\bigoplus\limits_{\substack{1 \leq r \leq a\\1 \leq i \leq \min{(r,d-a)}}} V_{r,i}^{a}\:.
\]
It then immediately follows from Theorem~\ref{thm:combinatorialdesign} that an $(m,a)$-ovoid of $\mathcal{P}$ is a combinatorial design with respect to the generators of $\mathcal{P}$.

The following lemma will be useful.
\begin{lemma}\label{lem:jprojection}
	Let $\bm{\chi}$ be the characteristic vector for an $(m,a)$-ovoid. Then there is a unique vector $\bm{u} \in {\left( V_{0,0}^{a}\right)}^{\perp}$
	for which
	\[ \bm{\chi} = \frac{m}{\qbin{d}{a}}\bm{j} + \bm{u}. \]
	where $\bm{u} \in {\left( V_{0,0}^{a}\right)}^{\perp}$.
\end{lemma}
\begin{proof}
	We can uniquely write $\bm{\chi} = s\bm{j} + \bm{u}$, with $\bm{u} \in {\left( V_{0,0}^{a}\right)}^{\perp}$.
	By the above remarks, $\bm{u} \in \ker{C_{a,d}}$ and so we have $C_{a,d}\bm{\chi} = s(C_{a,d}\bm{j}) = m\bm{j}$.
	Since each $d$-space contains $\qbin{d}{a}$ distinct $a$-spaces, we have $C_{a,d}\bm{j}=\qbin{d}{a}\bm{j}$. We can easily solve for $s$ to obtain the result.
\end{proof}

\begin{corollary}\label{cor:sizeovoid}
	An $(m,a)$-ovoid of $\mathcal{P}$ contains $m\prod_{i=1}^{a} \left(q^{d+e-i} +1 \right)$ elements.
\end{corollary}
\begin{proof}
	This follows immediately from the definition and Theorem~\ref{thm:subspacecount}. Alternatively, this can also be proven using Lemma~\ref{lem:jprojection}, taking the inner product with $\bm{j}$.
\end{proof}

\subsection{Regular \texorpdfstring{$(m,a)$}{(m,a)}-ovoids}

In this subsection we will impose an additional condition on $(m,a)$-ovoids, introducing regular ovoids, which we will need in Section~\ref{sec:newexample}.
Recall that the condition that a set $M$ with characteristic vector $\bm{\chi}$ is an $(m,a)$-ovoid is
\[
\bm{\chi} \in {\left( \bigoplus\limits_{1 \leq r \leq a} V_{r,0}^{a} \right)}^{\perp}
= V_{0,0}^{a}\ \oplus\bigoplus\limits_{\substack{1 \leq r \leq a\\1 \leq i \leq \min{(r,d-a)}}}V_{r,i}^{a}\:.
\]

We first look at the case $a=d-1$.
By the definition of the parameter $e$ we have that each $(d-1)$-space of $\mathcal{P}$ is contained in $q^{e}+1$ generators. Thus we can show the following theorem.
\begin{theorem}\label{thm:intersectionrelation}
	Let $\mathcal{P}$ be a 
	polar space of rank $d$ and parameter $e$, defined over $\F_{q}$. If $M$ is an $(m,d-1)$-ovoid of $\mathcal{P}$, then for any $\pi \in \Omega_{d-1}$ the number of subspaces $\rho \in M$, with $\rho \neq \pi$, and $\rho$ and $\pi$ contained in a common generator of $\mathcal{P}$ is given by
	\begin{align*}
		\begin{cases}
			(m-1)(q^{e}+1) & \mbox{ when } \pi \in M\:,     \\
			m(q^{e}+1)     & \mbox{ when } \pi \not\in M\:.
		\end{cases}
	\end{align*}
\end{theorem}
\begin{proof}
	For any $\pi \in \Omega_{d-1}$ there are $q^{e}+1$ generators through $\pi$.
	When $\pi \in M$, for each such generator $\omega$ there are $(m-1)$ other elements of $M$ contained in $\omega$.
	When $\pi \not\in M$, there are $m$ elements of $M$ contained in each of the $(q^{e}+1)$ generators through $\pi$.
	Furthermore there can at be at most one generator containing a fixed pair of $(d-1)$ spaces.
\end{proof}

\begin{remark}
	Using the relations introduced in Theorem~\ref{th:relations} we can write the previous theorem as follows: for any $\pi \in \Omega_{d-1}$:
	\[
	\left|\left\{\rho\in M:(\rho,\pi)\in R^{d-2,d}_{d-1,d-1}\right\}\right|=
	\begin{cases}
		(m-1)(q^{e}+1) & \text{ when } \pi \in M\:,     \\
		m(q^{e}+1)     & \text{ when } \pi \not\in M\:.
	\end{cases}
	\]
\end{remark}

Alternatively, we can prove the above result using the following corollary from Theorem~\ref{thm:eigenvaluesLie}.

\begin{corollary}
	For each $r$ with $1\leq r\leq d-1$ the eigenvalue of $C_{d-1,d-1}^{d-2,d}$ on the irreducible $V_{r,1}^{d-1}$ is $-(q^{e}+1)$.
\end{corollary}

We know if $M$ is an $(m,d-1)$-ovoid of $\mathcal{P}$ with characteristic vector $\bm{\chi}$, then $\bm{\chi}-\frac{m}{\qbin{d}{1}}\bm{j}$ is in the eigenspace of $C_{d-1,d-1}^{d-2,d}$ corresponding to the eigenvalue $-(q^{e}+1)$.

The relation $R_{a,a}^{a-1,a+1}$ corresponds to the graph of Lie type on the $a$-spaces, but these results show that we cannot use it, at least for the important case $a=d-1$, to discriminate between the $(m,a)$-ovoids.
Indeed, in Theorem~\ref{thm:intersectionrelation} we proved an intersection property with respect to the relation $R_{a,a}^{a-1,a+1}$ that holds for all $(m,d-1)$-ovoids.
For this reason, we turn to the relation $R_{a,a}^{a-1,a}$, which exists for all $a$. Imposing a a similar intersection condition leads to the definition of regular $(m,d-1)$-ovoids.
\begin{definition}\label{def:regularovoida}
	An $(m,a)$-ovoid $M$ of $\mathcal{P}$ will be called \emph{regular} if there are constants $c_{1}, c_{2}$ such that, for any $\pi \in \Omega_{a}$,
	\begin{align*}
		\left|\left\{\rho\in M:(\rho,\pi)\in R^{a-1,a}_{a,a}\right\}\right|=
		\begin{cases}
			c_{1} & \mbox{ when } \pi \in M,     \\
			c_{2} & \mbox{ when } \pi \not\in M.
		\end{cases}
	\end{align*}
\end{definition}

\begin{remark}
	It is clear that an $(m,a)$-ovoid whose characteristic vector belongs to $V_{0,0}^{a}\oplus V_{r,i}^{a}$ is regular.
\end{remark}

We can determine the possible values for $c_{1}$ and $c_{2}$. We immediately prove that the converse of the previous remark is also true.

\begin{theorem}\label{th:regularovoida}
	If $\bm{\chi}$ is the characteristic vector of a non-trivial regular $(m,a)$-ovoid $M$ of $\mathcal{P}$,
	then there is unique tuple $(r,i)$ with $1\leq r\leq a$ and $1\leq i\leq\min(r,d-a)$ such that $\bm{\chi} = \frac{m}{\qbin{d}{1}}\bm{j} + \bm{u}$ with $\bm{u}\in V^{a}_{r,i}$.
\end{theorem}
\begin{proof}
	We recall from Lemma~\ref{lem:jprojection} and the beginning of this subsection that
	\[
	\bm{\chi} = \frac{m}{\qbin{d}{1}}\bm{j} + \bm{u}\quad\text{ with }\quad \bm{u} \in V= \bigoplus\limits_{\substack{1 \leq r \leq a\\1 \leq i \leq \min{(r,d-a)}}}V_{r,i}^{a}\:.
	\]
	Moreover $\bm{u}\neq\bm{0}$ since $M$ is non-trivial. It follows from Definition~\ref{def:regularovoida} (also using its notation) that
	\[
	C_{a,a}^{a-1,a}\bm{\chi} = c_{1} \bm{\chi} + c_{2}(\bm{j} - \bm{\chi}) = (c_{1}-c_{2}) \bm{\chi} + c_{2}\bm{j}\:.
	\]
	So, we find
	\begin{align*}
		C_{a,a}^{a-1,a}\bm{u} & = \left(c_{1}-c_{2}\right)\bm{u} + \left(c_{2}+\frac{m}{\qbin{d}{1}}\left(c_{1}-c_{2}-\mu_{d,a,a-1,a,0,0}\right)\right)\bm{j}\:,
	\end{align*}
	since $\bm{j}\in V^{a}_{0,0}$. Since $\bm{u} \in V$ we know that the coefficient of $\bm{j}$ equals 0, and that $\bm{u}$ is an eigenvector of $C_{a,a}^{a-1,a}$.
	From Theorem~\ref{th:uniqueeigenvalues} it follows that $\bm{u}\in V^{a}_{r,i}$ for a unique $r$ and $i$.
\end{proof}

\begin{remark}
	We choose to call these $(m,a)$-ovoids regular because a set of $a$-spaces whose characteristic vector satisfies the above mentioned condition
	will be a \emph{regular set} with respect to each relation of the association scheme on $a$-spaces,
	that is, the condition of Definition~\ref{def:regularovoida}
	will in fact hold for every relation $R_{a,a}^{s,k}$.
\end{remark}

Given the previous theorem, the following definition is well-defined.

\begin{definition}\label{def:type}
	Let $\mathcal{P}$ be a polar space of rank $d$. If the characteristic vector of a non-trivial regular $(m,a)$-ovoid $M$ is contained in $V^{a}_{0,0}\oplus V^{a}_{r,i}$, then it has \emph{type} $(r,i)$.
	In cases $a=d-1$ or $a=1$, it is necessarily of type $(r,1)$ and we will shorten this to type $r$.
\end{definition}

\begin{example}\label{ex:embeddedovoid}
	Let $\mathcal{P}$ be a 
	polar space of rank $d$ and parameter $e\leq1$, defined over $\F_{q}$. If $\mathcal{P}$ is a quadric or a Hermitian variety, there is an embedded 
	polar $\mathcal{P}'$ space of rank $d-1$ and parameter $e+1$. It arises as the intersection of $\mathcal{P}$ with a non-tangent hyperplane of the ambient projective space.
	The set of $a$-spaces of $\mathcal{P}'$ is a regular $\left(\qbin{d-1}{a},a\right)$-ovoid of type $(1,1)$ of $\mathcal{P}$, for any $a$ with $1\leq a\leq d-1$. This follows from~\cite[Theorem 4.4.3]{phdvanhove}.
	\par We present all embeddings that can be used in this construction.
	\begin{itemize}
		\item There are parabolic quadrics $\mathcal{Q}(2d,q)$ embedded in a hyperbolic quadric $\mathcal{P}=\mathcal{Q}^{+}(2d+1,q)$.
		\item There are elliptic quadrics $\mathcal{Q}^{-}(2d-1,q)$ embedded in the parabolic quadric $\mathcal{P}=\mathcal{Q}(2d,q)$.
		\item There are Hermitian polar spaces $\mathcal{H}(2d,q)$ embedded in the Hermitian polar space $\mathcal{P}=\mathcal{H}(2d+1,q)$, $q$ a square.
	\end{itemize}
\end{example}

\begin{remark}\label{rem:ovoida1}
	An $(m,1)$-ovoid, which is an $m$-ovoid in the traditional sense, i.e.~a point set of $\mathcal{P}$ such that any generator of $\mathcal{P}$ contains $m$ points of it, is always a regular $(m,1)$-ovoid of type 1.
	Indeed, for $a=1$ we have that
	\[
	\bigoplus\limits_{\substack{1 \leq r \leq a\\1 \leq i \leq \min{(r,d-a)}}}V_{r,i}^{a}=V_{1,1}^{1}.
	\]
\end{remark}

We now look at the particular case $a=d-1$, and determine the intersection numbers of regular $(m,d-1)$-ovoids with respect to all relations. We use the eigenvalues calculated above.

\begin{lemma}\label{lem:intersectionnumbersgeneralk}
	Let $M$ be a non-trivial regular $(m,d-1)$-ovoid $M$ of type $r$ in $\mathcal{P}$.
	Then, for any $\pi \in \Omega_{d-1}$, the number of subspaces $\rho \in M$, with $\dim(\rho\cap\pi)=s$ and $\dim(\langle \rho, \rho^\perp\cap \pi\rangle)=k$, with $0\leq s\leq d-1$ and $d-1\leq k\leq d$ is given by
	\begin{align*}
		\begin{cases}
			\mu_{d,d-1,s,k,r,1}+\frac{m}{\qbin{d}{1}}(\mu_{d,d-1,s,k,0,0}-\mu_{d,d-1,s,k,r,1}) & \quad \text{ if } \pi \in M\:,   \\
			\frac{m}{\qbin{d}{1}}(\mu_{d,d-1,s,k,0,0}-\mu_{d,d-1,s,k,r,1})                     & \quad\text{ if } \pi \notin M\:.
		\end{cases}
	\end{align*}
\end{lemma}
\begin{proof}
	Let $\bm{\chi}$ be the characteristic vector of the non-trivial regular $(m,d-1)$-ovoid $M$ of type $r$. From Definition~\ref{def:type} it follows that
	\[
	\bm{\chi} = \frac{m}{\qbin{d}{1}}\bm{j} + \bm{u}\quad\text{ with }\quad \bm{u} \in V_{r,1}^{d-1}\:.
	\]
	We find that
	\begin{align*}
		C_{d-1,d-1}^{s,k}\bm{\chi} & =\frac{m}{\qbin{d}{1}}C_{d-1,d-1}^{s,k}\bm{j} + C_{d-1,d-1}^{s,k}\bm{u}=\frac{m}{\qbin{d}{1}}\mu_{d,d-1,s,k,0,0}\bm{j} + \mu_{d,d-1,s,k,r,1}\bm{u} \\
		& =\frac{m}{\qbin{d}{1}}\mu_{d,d-1,s,k,0,0}\bm{j} + \mu_{d,d-1,s,k,r,1}\left(\bm{\chi} - \frac{m}{\qbin{d}{1}}\bm{j}\right)                          \\
		& =\frac{m}{\qbin{d}{1}}(\mu_{d,d-1,s,k,0,0}-\mu_{d,d-1,s,k,r,1})((\bm{j}-\bm{\chi})+\bm{\chi}) + \mu_{d,d-1,s,k,r,1}\bm{\chi}                       \\
		& =\frac{m}{\qbin{d}{1}}(\mu_{d,d-1,s,k,0,0}-\mu_{d,d-1,s,k,r,1})(\bm{j}-\bm{\chi}) +                                                                \\&\qquad\left(\mu_{d,d-1,s,k,r,1}+\frac{m}{\qbin{d}{1}}(\mu_{d,d-1,s,k,0,0}-\mu_{d,d-1,s,k,r,1})\right)\bm{\chi}\:.
	\end{align*}
	This shows that for any $\pi \in \Omega_{d-1}$, the number of subspaces $\rho \in M$ which are in relation $R^{s,k}_{d-1,d-1}$ only depends on whether $\pi$ itself belongs to $M$. The numbers are the coefficients of $\bm{j}-\bm{\chi}$ and $\bm{\chi}$.
\end{proof}

Note that the actual numbers and can be computed using Theorem~\ref{th:eigenval} (which provides $\mu_{d,d-1,s,k,r,1}$) and Corollary~\ref{countingBCN} (which provides $\mu_{d,d-1,s,k,0,0}$). We present these numbers for the regular ovoids of type 1, i.e.~the case $r=1$.

\begin{corollary}\label{cor:intersectionnumbersgeneralk}
	Let $M$ be a non-trivial regular $(m,d-1)$-ovoid $M$ of type $1$ in $\mathcal{P}$. Then, for any $\pi \in \Omega_{d-1}$,
	the number of subspaces $\rho \in M$, with $\dim(\rho\cap\pi)=s$ and $\dim(\langle \rho, \rho^\perp\cap \pi\rangle)=k$, with $0\leq s\leq d-1$ and $d-1\leq k\leq d$, equals
	\begin{align*}
		\begin{cases}
			q^{e(d-s-1)+\binom{d-s}{2}-1}\left((m(q-1)+1)\qbin{d-2}{s}+q^{d-s}\qbin{d-2}{s-1}\right) & \quad \text{ if } \pi \in M\:,   \\
			m \left(q-1\right)q^{e(d-s-1)+\binom{d-s}{2}-1}\qbin{d-2}{s}                             & \quad\text{ if } \pi \notin M\:.
		\end{cases}
	\end{align*}
	if $k=d-1$, and equals
	\begin{align*}
		\begin{cases}
			(m-1)(q^e+1)q^{e(d-s-2)+\binom{d-s}{2}-1}\qbin{d-2}{s}          & \quad \text{ if } \pi \in M\:,   \\
			m \left(q^e+1\right) q^{e(d-s-2)+\binom{d-s}{2}-1}\qbin{d-2}{s} & \quad\text{ if } \pi \notin M\:.
		\end{cases}
	\end{align*}
	if $k=d$.
\end{corollary}
\begin{proof}
	Note that it follows from Corollary~\ref{countingBCN} that
	\begin{align*}
		\mu_{d,d-1,s,k,0,0}=
		\begin{cases}
			q^{e(d-s-1)+\binom{d-s+1}{2}-1}\qbin{d-1}{s}                                & \quad k=d-1\:, \\
			q^{e(d-s-2)+\binom{d-s}{2}}\qbin{d-1}{s}\qbin{d-1-s}{1}\left(q^{e}+1\right) & \quad k=d\:,
		\end{cases}
	\end{align*}
	and from Theorem~\ref{th:eigenval} that
	\begin{align*}
		\mu_{d,d-1,s,k,1,1}=
		\begin{cases}
			q^{e(d-s-1)+\binom{d-s}{2}-1}\left(\qbin{d-2}{s}+\qbin{d-2}{s-1}q^{d-s}\right) & \quad k=d-1\:, \\
			-\qbin{d-2}{s}q^{e(d-s-2)+\binom{d-s}{2}-1}(q^{e}+1)                           & \quad k=d\:.
		\end{cases}
	\end{align*}
	The result then follows from Lemma~\ref{lem:intersectionnumbersgeneralk}
\end{proof}


\subsection{Classification of the regular \texorpdfstring{$(1,d-1)$}{(1,d-1)}-ovoids of type 1}

The regular $(m,d-1)$-ovoids of type 1 are of importance in Section~\ref{sec:newexample} and for that reason we investigate them in more detail in this subsection.
In particular, we focus on the case $m=1$. We will prove that all $(1,d-1)$-ovoids of type 1 essentially arise from Example~\ref{ex:embeddedovoid}.
The proof for polar spaces of rank 3 and of rank at least 4 is slightly different.

\begin{theorem}\label{th:classification1ovoids}
	Every $(1,d-1)$-ovoid of type 1 in a polar space of rank $d\geq3$ and with parameter $e\leq1$ is the set of generators of an embedded polar space of rank $d-1$ with parameter $e+1$,
	hence corresponds to Example~\ref{ex:embeddedovoid} for $a=d-1$.
\end{theorem}
\begin{proof}
	Let Let $M$ be a regular $(1,d-1)$-ovoid of type 1 in $\P$.
	Then we know by Corollary~\ref{cor:intersectionnumbersgeneralk} that for any $(d-1)$-space $\pi$ in $\P$ and for any $s=0,\dots,d-1$ we have
	\begin{align}
		\left|\left\{\rho\in M:(\rho,\pi)\in R^{s,d-1}_{d-1,d-1}\right\}\right| & =
		\begin{cases}
			q^{e(d-s-1)+\binom{d-s}{2}}\qbin{d-1}{s}                   & \text{ if } \pi \in M\:,     \\
			\left(q-1\right)q^{e(d-s-1)+\binom{d-s}{2}-1}\qbin{d-2}{s} & \text{ if } \pi \not\in M\:,
		\end{cases}\label{eq:regularovoidclass1} \\
		\left|\left\{\rho\in M:(\rho,\pi)\in R^{s,d}_{d-1,d-1}\right\}\right|   & =
		\begin{cases}
			0                                                            & \text{ if } \pi \in M\:,     \\ 
			\left(q^e+1\right)q^{e(d-s-2)+\binom{d-s}{2}-1}\qbin{d-2}{s} & \text{ if } \pi \not\in M\:.
		\end{cases}\label{eq:regularovoidclass2}
	\end{align}
	
	We know that $|M|=\prod_{i=1}^{d-1} \left(q^{e+i} +1 \right)$ by Corollary~\ref{cor:sizeovoid}. Let $\pi$ be a $(d-1)$-space in $M$. We will show by induction that for any $j=0,\dots,d-1$,
	through any $(d-1-j)$-space $\pi_{d-j-1}$ in $\pi$ there are precisely 
	$q^{(e+1)j+\binom{j}{2}}$ distinct $(d-1)$-spaces of $M$ whose intersection with $\pi$ equals $\pi_{d-j-1}$, and that for each generator through $\pi_{d-j-1}$ the unique element of $M$ in it contains $\pi_{d-j-1}$.
	For $j=0$ the statement is immediate.
	\par We now assume the statement is true for $0,\dots,j-1$ and we show it for $j$. Let $\pi_{d-j-1}$ be a $(d-1-j)$-space in $\pi$.
	We know that that the tangent space $T_{\pi_{d-j-1}}(\P)$ meets $\P$ in a cone with vertex $\pi_{d-j-1}$ and basis a polar space $\P_{j}$ of rank $j+1$ and of the same type as $\P$.
	The elements of $M$ through $\pi_{d-j-1}$ correspond to $j$-spaces of $\P_{j}$, and each generator (i.e.~$(j+1)$-space) of $\P_{j}$ contains at most one of these $j$-spaces, since $M$ is a $(1,d-1)$-ovoid.
	So, the number of $(d-1)$-spaces of $M$ through $\pi_{d-j-1}$ is at most $\prod_{i=1}^{j} \left(q^{e+i} +1 \right)$, which is the size of a $(1,j)$-ovoid in $\P_{j}$ by Corollary~\ref{cor:sizeovoid}.
	For any $k=1,\dots,j$, there are $\qbin{j}{k}$ distinct $(d-j-1+k)$-spaces through $\pi_{d-j-1}$ in $\pi$, and by the induction hypothesis,
	through each $(d-j-1+k)$-space there are $q^{(e+1)(j-k)+\binom{j-k}{2}}$ elements of $M$ whose intersection with $\pi$ equals this $k$-space; so there are
	\begin{align*}
		\sum_{k=1}^{j}\qbin{j}{k}q^{(e+1)(j-k)+\binom{j-k}{2}} & = \sum_{k=0}^{j}\qbin{j}{k}q^{(e+1)k+\binom{k}{2}}-q^{(e+1)j+\binom{j}{2}} \\
		& = \prod^{j-1}_{k=0}\left(1+q^{e+1+k}\right)-q^{(e+1)j+\binom{j}{2}}        \\
		& = \prod^{j}_{i=1}\left(1+q^{e+i}\right)-q^{(e+1)j+\binom{j}{2}}
	\end{align*}
	elements of $M$ through $\pi_{d-j-1}$ meeting $\pi$ in a subspace of dimension at least $d-j$, where we used Lemma~\ref{lem:gaussianbinomialidentity} in the second step.
	Hence, the number of $(d-1)$-spaces of $M$ through $\pi_{d-j-1}$ whose intersection with $\pi$ is precisely $\pi_{d-j-1}$, is at most $q^{(e+1)j+\binom{j}{2}}$.
	Since $\pi$ contains $\qbin{d-1}{d-j-1}$ distinct $(d-j-1)$-spaces, and $\pi$ meets $q^{ej+\binom{j+1}{2}}\qbin{d-1}{d-j-1}+0=q^{(e+1)j+\binom{j}{2}}\qbin{d-1}{d-j-1}$ elements of $M$
	in precisely a $(d-j-1)$-space by~\eqref{eq:regularovoidclass1} and~\eqref{eq:regularovoidclass2},
	each $(d-j-1)$-space of $\pi$ must be on exactly $q^{(e+1)j+\binom{j}{2}}$ elements of $M$ meeting $\pi$ in precisely this $(d-j-1)$-space,
	hence on $\prod^{j}_{i=1}\left(1+q^{e+i}\right)$ distinct $(d-1)$-spaces of $M$ meeting $\pi$ in a subspace of dimension at least $d-j-1$.
	In other words, a $(d-j-1)$-space is on zero or on $\prod^{j}_{i=1}\left(1+q^{e+i}\right)$ distinct $(d-1)$-spaces of $M$, and is called poor or rich, respectively.
	In particular the elements of $M$ are precisely the rich $d-1$-spaces.
	Moreover, each generator through a rich $(d-j-1)$-space $\sigma$ contains precisely one $(d-1)$-space of $M$ through $\sigma$, since the corresponding $j$-spaces in $\P_{j}$ form a $(1,j)$-ovoid.
	In other words, for any generator through $\sigma$, the unique $(d-1)$-space of $M$ in it, goes through $\sigma$.
	\medskip\par
	\medskip\par Let $\ell$ be a poor line, and assume that there are two rich points on it. Then for any generator through $\ell$, the elements of $M$ in it should go through both of the rich points, and thus also through the line, a contradiction. So, a line containing two rich points must be a rich line.
	\par Now, consider a rich point $P$ and a rich line $\ell$ not through $P$. If $P$ and $\ell$ generate a plane $\pi$ of $\mathcal{P}$, then for every generator through $\pi$,
	the unique element of $M$ in it should contain both $P$ and $\ell$ and hence the plane they span. So, in this case, for every point $Q$ on $\ell$ the line through $P$ and $Q$ is rich. Note that this case does not occur if $d=3$.
	\par If $P$ and $\ell$ do not generate a plane in $\mathcal{P}$, then there is a unique line $\ell'$ on $\mathcal{P}$ through $P$ meeting $\ell$.
	The intersection point $P'$ of $\ell$ and $\ell'$ is a rich point since it lies on the rich line $\ell$. As $\ell'$ contains two rich points, it is a rich line. So, in this case there is a unique rich line $\ell'$ through $P$, and meeting $\ell$ in a point.
	\par The geometry of the rich points and rich lines thus fulfills the `one or all' axiom and is thus a Shult space, and it is non-degenerate, since a rich point is not on all elements of $M$.
	By the Buekenhout-Lefèvre result~\cite{bl} (see also~\cite[Theorem 26.3.29]{Hirschfeld3}) a non-degenerate Shult space embedded in a finite projective space is a finite classical polar space.
	So, $M$ is the set of generators of a polar space of rank $d-1$ and with parameter $e+1$ embedded in $\P$.
\end{proof}


\section{A new Cameron-Liebler set in polar spaces with \texorpdfstring{$e\leq1$}{e=<1}}\label{sec:newexample}

\begin{example}\label{ex:CLps}
	Consider a polar space $\mathcal{P}$ of rank $d+2\geq4$ with parameter $e\leq1$ in $\PG(n+3,q)$,
	and let $\alpha\cong\PG(n,q)$ be an $(n+1)$-dimensional subspace that intersects $\mathcal{P}$ in a polar space $\mathcal{P}'$ of rank $d$ with parameter $e+1$.
	So, if $n=2d$, then $\mathcal{P}$ is a hyperbolic quadric $\cQ^{+}(2d+3,q)$ or a Hermitian polar space $\cH(2d+3,q)$, with $\mathcal{P}'$ a parabolic quadric $\cQ(2d,q)$ or a Hermitian polar space $\cH(2d,q)$, respectively.
	In the Hermitian case we assume $q$ to be a square. If $n=2d+1$, then $\mathcal{P}$ is a parabolic quadric $\cQ(2d+4,q)$, with $\mathcal{P}'$ an elliptic quadric $\cQ^{-}(2d+1,q)$.
	\par Let $\mathcal{M}$ be a regular $(m,d-1)$-ovoid of type 1 in $\mathcal{P}'$, see Definitions~\ref{def:maovoid},~\ref{def:regularovoida} and~\ref{def:type}.
	So, this is a set of $(d-1)$-spaces in $\mathcal{P}'$ (which has parameter $e+1$) such that any generator of $\mathcal{P}'$ contains precisely $m$ elements of $\mathcal{M}$, and such that for a given $(d-1)$-space $\sigma\in\mathcal{P}'$,
	the number of $(d-1)$-spaces $\sigma_0$ of $\mathcal{M}$ meeting $\sigma$ in a $(d-2)$-space (space of codimension $1$), and with $\langle \sigma, \sigma_0\rangle$ not a generator of $\mathcal{P}'$ equals
	\begin{align*}
		\begin{cases}
			mq^{e+1}(q-1) & \text{ if } \sigma \notin \mathcal{M}\:, \\
			(m-1)q^{e+1}(q-1)+ q^{e+2}\qbin{d-1}{1}
			& \text{ if } \sigma \in \mathcal{M}\:.
		\end{cases}
	\end{align*}
	\par Let $\L$ be the family of generators of $\mathcal{P}$ that meet $\alpha$ in a subspace of dimension $d-1$ and such that this intersection is an element of $\mathcal{M}$.
\end{example}

\begin{theorem}\label{CL9}
	Using the notation from Example~\ref{ex:CLps}, the set $\L$ is a Cameron-Liebler set of generators in $\mathcal{P}$ with parameter $mq^{e+1}(q-1)$.
\end{theorem}
\begin{proof}
	We will apply Theorem~\ref{stellinggg} for $i=1$. So, it will be sufficient to prove that for every generator $\pi$ of $\P$, the number of elements of $\L$ meeting $\pi$ in a $(d+1)$-space (space of codimension 1) is equal to
	\begin{align*}
		\begin{cases}
			mq^{e+1}(q-1)-1+q^{e}\frac{q^{d+1}-1}{q-1} \quad & \text{if } \pi \in \L\:,    \\
			mq^{e+1}(q-1) \qquad                             & \text{if } \pi \notin \L\:.
		\end{cases}
	\end{align*}
	Recall that the generators of $\mathcal{P}$ are $(d+2)$-spaces.
	
	\paragraph{Case 1: $\pi\in \L$}\mbox{}\medskip
	
	The generator $\pi$ meets $\alpha$ in a $(d-1)$-space $\beta\in\mathcal{M}$. We first count the number of generators in $\L$, containing $\beta$ and meeting $\pi$ in a $(d+1)$-space.
	These generators and also $\pi$ are contained in $\mathcal{P}\cap T_{\beta}(\mathcal{P})$, which is a cone with vertex $\beta$ and base a polar space $\mathcal{P}_{3}$ of rank 3 and of the same type as $\mathcal{P}$.
	The tangent space $T_{\beta}(\mathcal{P}')$ is contained in $T_{\beta}(\mathcal{P})$, and has as base a curve $C$ of the same type as $\mathcal{P}'$, so a $\cQ(2,q)$, a $\cH(2,q)$ or a $\cQ^{-}(3,q)$.
	We may assume that $C$ is contained in $\mathcal{P}_{3}$. The generator $\pi$ corresponds to a plane $\pi'$, disjoint from $C$, and the generators of $\L$ through $\beta$ meeting $\pi$ in a $(d+1)$-space correspond to planes in $\mathcal{P}_{3}$, meeting $\pi'$ in a line and disjoint from $C$.
	The number of planes meeting $\pi'$ in a line is $(q^2+q+1)q^{e}$. Through a point of $C$ there is one plane in $\mathcal{P}_{3}$ meeting $\pi'$ in a line.
	So, the number of planes in $\mathcal{P}_{3}$, meeting $\pi'$ in a line and disjoint from $C$ equals $(q^2+q+1)q^{e}-|C|=(q^2+q+1)q^{e}-(q^{e+1}+1)=(q^{2}+1)q^{e}-1$.
	\par Now we count the number of generators in $\L$, not containing $\beta$ and meeting $\pi$ in a $(d+1)$-space. Such a generator meets $\alpha$ in a $(d-1)$-space belonging to $\mathcal{M}$, which meets $\beta$ in a $(d-2)$-space.
	We know, by the definition of $\mathcal{M}$, that there are $(m-1)q^{e+1}(q-1)+ q^{e+2}\qbin{d-1}{1}$ elements $\beta'$ of $\mathcal{M}$ meeting $\beta$ in a $(d-2)$-space,
	and such that $\langle \beta, \beta'\rangle$ is not a generator of $\mathcal{P}'$.
	For each such $\beta'$, there is precisely one generator though it, which meets $\pi$ in a $(d+1)$-space. Since $\langle \beta, \beta'\rangle\not\subset\mathcal{P'}$, this generator does not meet $\alpha$ in a $d$-space, so it belongs to $\L$.
	\par Also, for every $\beta'\in \mathcal{M}$ with $\langle \beta, \beta'\rangle$ a generator of $\mathcal{P}'$ and meeting $\beta$ in a $(d-2)$-space,
	we know that there is precisely one generator through $\beta'$, meeting $\pi$ in a $(d+1)$-space, but this generator contains $\langle \beta, \beta'\rangle$, and hence, is not contained in $\L$ since it meets $\alpha$ in a $d$-space.
	
	We conclude that there are
	\begin{multline*}
		(q^{2}+1)q^{e}-1+(m-1)q^{e+1}(q-1)+ q^{e+2}\qbin{d-1}{1} 
		= mq^{e+1}(q-1)-1+q^{e}\qbin{d+1}{1}
	\end{multline*}
	generators in $\mathcal{L}$ meeting $\pi$ in a $(d+1)$-space.
	
	\paragraph{Case 2: $\pi\notin \L$ and $\pi$ contains a $(d-1)$-space of $\mathcal{M}$}\mbox{}\medskip
	
	As we know that $\pi\notin\L$, we know that $\pi$ meets $\mathcal{P}'$ in a $d$-space (i.e.~a generator of $\mathcal{P}'$), which contains $m$ distinct $(d-1)$-spaces belonging to $\mathcal{M}$, say $\beta_1, \dots, \beta_m$.
	A generator of $\L$ meeting $\pi$ in a $(d+1)$-space should contain $\beta_{i}$ for some $i$. Given a $(d-1)$-space $\beta_i\in\mathcal{M}$, we count the number of elements $\gamma$ through it meeting $\pi$ in a $(d+1)$-space.
	For this, we consider the tangent space $T_{\beta_i}(\mathcal{P})$ of $\beta_i$ to $\mathcal{P}$, which is a cone with vertex $\beta_i$ and base a polar space $\mathcal{P}_{3}$ of rank 3 and of the same type as $\mathcal{P}$.
	The tangent space $T_{\beta_i}(\mathcal{P}')$ is contained in $T_{\beta_i}(\mathcal{P})$, and has as base a curve $C$ of the same type as $\mathcal{P}'$, so a $\cQ(2,q)$, a $\cH(2,q)$ or a $\cQ^{-}(3,q)$.
	We may assume that $C$ is contained in $\mathcal{P}_{3}$. The generator $\pi$ corresponds to a plane $\pi'$ in $\mathcal{P}_{3}$, that contains a point of $C$.
	A generator $\gamma\in\L$ through $\beta_{i}$ meeting $\pi$ in a $(d+1)$-space corresponds to a plane $\gamma'$ in $\mathcal{P}_{3}$ that meets $\pi$ in a line and is disjoint from $C$.
	So, we count the number of planes in $\mathcal{P}_{3}$ meeting $\pi'$ in a line, and disjoint to $C$. There are $q^2$ lines in $\pi$ disjoint to the point $C\cap\pi'$, and each of them is contained in $q^{e}$ planes, different from $\pi'$.
	Through each of the $q^{e+1}$ points of $C\setminus\pi'$ there is precisely one plane in $\mathcal{P}_{3}$ meeting $\pi'$ in a line. Hence, through each $\beta_{i}$ there are $q^{e+2}-q^{e+1}$ generators in $\L$ meeting $\pi$ in a $(d+1)$-space.
	So, in total there are $mq^{e+1}(q-1)$ generators in $\cL$, meeting $\pi$ in a $(d+1)$-space.
	
	\paragraph{Case 3: $\pi\notin \L$ and $\pi$ does not contain an $(d-1)$-space of $\mathcal{M}$}\mbox{}\medskip
	
	Let $\beta$ be the $(d-1)$-space $\pi\cap \mathcal{P}'$.
	Any generator in $\L$ meeting $\pi$ in a $(d+1)$-space should meet $\beta$ in a subspace of dimension at least $(d-2)$.
	For every $(d-1)$-space $\mu\in \mathcal{M}$ meeting $\beta$ in a $(d-2)$-space, and with $\langle \beta, \mu\rangle\nsubseteq \mathcal{P}'$, we know that there is precisely one generator through $\mu$ meeting $\pi$ in a $(d+1)$-space.
	The number of $(d-1)$-spaces $\mu$ follows from the definition of a regular $(m,d-1)$-ovoid, and is equal to $mq^{e+1}(q-1)$.
	For every $(d-1)$-space $\mu\in \mathcal{M}$ meeting $\beta$ in a $(d-2)$-space, and with $\langle \beta, \mu\rangle\subseteq \mathcal{P}'$, we know that the only generator through $\mu$ meeting $\pi$ in a $(d+1)$-space necessarily contains $\langle\beta,\mu\rangle$.
	Hence, it meets $\alpha$ in a subspace of dimension $d$ and is thus not a generator in $\L$. So, we find $mq^{e+1}(q-1)$ elements of $\L$ meeting $\pi$ in a $(d+1)$-space.
\end{proof}

\begin{remark}\label{rem:newconstructiontrivial}
	We use the notation of Example~\ref{ex:CLps}. Assume that $\mathcal{M}$ is the set of generators ($(d-1)$-spaces) of an embedded polar space $\mathcal{P}''$ of rank $d-1$ and parameter $e+2$ in $\mathcal{P}'$.
	In particular we know this implies that $e=0$, and that $\mathcal{P}$ is a $\mathcal{Q}^{+}(2d+3,q)$, $\mathcal{P}'$ is a $\mathcal{Q}(2d,q)$ and $\mathcal{P}''$ is a $\mathcal{Q}^{-}(2d-1,q)$.
	By Theorem~\ref{th:classification1ovoids} we know that any regular $(1,d-1)$-ovoid in $\mathcal{P}'$ is such an embedded polar space if $d\geq 3$. The constructed Cameron-Liebler set $\mathcal{L}$ has parameter $q(q-1)$.
	\par Now, let $\alpha'$ be the hyperplane of $\alpha$, i.e.~an $n$-space, wherein $\mathcal{P}''$ is contained.
	The image $\alpha'^{\perp}$ of $\alpha'$ under the polarity $\perp$ associated with $\mathcal{P}$ is a 4-space intersecting $\alpha'$ trivially, which contains the 3-space $\alpha^{\perp}$.
	Moreover, $\alpha'^{\perp}\cap\mathcal{P}$ is an embedded $Q_{3}\cong\mathcal{Q}^{-}(3,q)$, and $\alpha^{\perp}\cap\mathcal{P}$ is a $Q_2\cong\mathcal{Q}(2,q)$ therein.
	\par We note that for a point $P\in Q_{3}\setminus Q_2$ the point-pencil with vertex $P$ is contained in $\mathcal{L}$.
	Since the points of $Q_{3}\setminus Q_2$ form a partial ovoid, we find that in this case the (degree one) Cameron-Liebler set $\mathcal{L}$ is the disjoint union of $q(q-1)$ point-pencils, and thus is trivial.
\end{remark}

Given the previous remark, it is important to check whether in general the (degree one) Cameron-Liebler sets constructed in Example~\ref{ex:CLps} do not contain trivial degree one Cameron-Liebler sets. We have the following results.

\begin{theorem}
	Let $\mathcal{P}$, $\mathcal{P}'$, $\mathcal{M}$ and $\mathcal{L}$ be as in Example~\ref{ex:CLps}.
	If $\mathcal{M}$ does not contain the set of generators of an embedded polar space of degree $d-1$ and parameter $e+2$ inside $\mathcal{P}'$, then the Cameron-Liebler set $\L$ of generators of $\P$ from Example~\ref{ex:CLps}  does not contain a point-pencil.
\end{theorem}
\begin{proof}
	Suppose that $\L$ contains a point-pencil $\L_P$ with vertex the point $P$. Then we know that $P\notin \P'$, as otherwise $\L$ would contain all generators of $\P$ through the point $P\in \P'$, and hence, as well the generators of $\P$ through $P$ meeting $\P'$ in a $d$-space, which contradicts the definition of $\L$.
	\par Now consider the intersection $T_P(\P)\cap \P$. Since $\L$ contains the point-pencil $\L_P$ by assumption, we know that all generators in $T_P(\P)\cap \P$ are contained in $\L$.
	If $\alpha\cap T_P(\P)$ contains a $d$-space of $\P$, then $\L_P\subset \L$ contains generators meeting $\alpha$ in a $d$-space, a contradiction.
	So, $\alpha\cap T_P(\P)$ should be an $n$-space that meets $\P$ in a non-singular polar space $\P''$ with rank $d-1$, and hence with parameter $e+2$.
	Since all of the generators of $T_{P}(\P)$ are contained in $\L$, they must meet $\alpha$ in precisely a $d-1$ space of $\mathcal{M}$.
	Therefore we find that $\mathcal{M}$ must contain the set of generators of $\P''$, contradicting the assumptions in the statement of the theorem.
\end{proof}

\begin{theorem}
	Let $\mathcal{P}$, $\mathcal{P}'$, $\mathcal{M}$ and $\mathcal{L}$ be as in Example~\ref{ex:CLps}.
	If $\mathcal{M}$ is not an embedded polar space of degree $d-1$ and parameter $e+2$ inside $\mathcal{P}'$, then the Cameron-Liebler set $\L$ of generators of $\P$ from Example~\ref{ex:CLps}  does not contain an embedded polar space as defined in Example in~\ref{ex:embedded}.
\end{theorem}
\begin{proof}
	If the polar space $\P$ has parameter $e<1$, then it does not contain an embedded polar space of the same rank.
	Hence, we only have to investigate the case $e=1$, in which $\P \equiv Q(2d+4,q)$ and $\P' \equiv Q^-(2d+1,q)$.
	Suppose that $\L$ contains an embedded polar space $\Tilde{\P}\equiv Q^+(2d+3,q)$ from Example~\ref{ex:embedded}.
	Then $\Tilde{\P}$ meets $\alpha$ in a (possibly degenerate) polar space of rank at least $d$, and hence $\Tilde{\P}$ contains generators meeting $\P'$ in a $d$-space.
	As $\Tilde{\P}\subset \L$, we have that $\L$ also contains the elements of $\Tilde{\P}$ meeting $\P'$ in a $d$-space, which cannot happen due to the definition of $\L$.
\end{proof}

\medskip

\par We conclude from the previous two theorems that the construction presented in Example~\ref{ex:CLps} is `new' in almost all cases.
Unless the regular $(m,d-1)$-ovoid used in the construction contains a regular $(1,d-1)$-ovoid that arises from an embedded polar space, we get a Cameron-Liebler set that was not known before, nor decomposes in Cameron-Liebler sets of which one or more are trivial.
\par Now we look at a few particular instances  of Example~\ref{ex:CLps}. The first one is clearly trivial, but for a parabolic quadric its complement has a nice description.

\begin{example}
	Using the notation of Example~\ref{ex:CLps}, the set of all $(d-1)$-spaces of $\P'$ is a trivial regular $\left(\qbin{d}{1},d-1\right)$-ovoid $\mathcal{M}$.
	As it is trivial, its type is not well-defined, but we can still apply the construction and Theorem~\ref{CL9}. Then $\mathcal{L}$ is the set of all generators of $\P$ whose intersection with $\alpha$ is a $(d-1)$-space in $\P'$.
	Its complement, $\overline{\mathcal{L}}$ is the set of all generators of $\P$ whose intersection with $\alpha$ is a $d$-space in $\P'$ (i.e.~a generator of $\P'$).
	Then both $\mathcal{L}$ and $\overline{\mathcal{L}}$ are (degree one) Cameron-Liebler sets of $\P$, with parameters $q^{e+1}(q^{d}-1)$ and $q^{e+1}+1$, respectively.
	\par We look in detail at $\overline{\mathcal{L}}$ because of its small parameter. If $\P=\mathcal{Q}^{+}(2d+3,q)$ or $\P=\mathcal{H}(2d+3,q)$, then $\alpha^\perp\cap \P$ is a $Q(2,q)$ or a $H(2,q)$, respectively, and $\overline{\mathcal{L}}$ is the disjoint union of all the point-pencils with vertex a point in $\alpha^\perp\cap \P$.
	If $\P=\mathcal{Q}(2d+4,q)$, then $\overline{\mathcal{L}}$ has parameter $q^{2}+1$. We know that $\alpha^\perp$ is a plane and that $C=\alpha^\perp\cap \P$ is a $Q(2,q)$.
	A hyperplane of $\PG(2d+4,q)$ through $\alpha$ meets $\P$ in a $\cQ^{-}(2d+3,q)$, a $\cQ^{+}(2d+3,q)$ or a cone with vertex a point and base a $\cQ(2d+2,q)$, if it meets $\alpha^\perp$ in a secant line to $C$, a line skew to $C$ or a tangent line to $C$, respectively.
	Each generator of $\P$ whose intersection with $\alpha$ is a $d$-space, spans a hyperplane with $\alpha$, and this hyperplane cannot meet $\P$ in an elliptic quadric $\cQ^{-}(2d+3,q)$.
	Hence, all the elements of $\overline{\mathcal{L}}$ are contained in a hyperplane that meets $\P$ in a hyperbolic quadric $\cQ^{+}(2d+3,q)$ or in a cone.
	Hence $\P$ is the disjoint union of $q+1$ point-pencils and $\frac{q(q-1)}{2}$ embedded hyperbolic quadrics, so it is a trivial example.
\end{example}

Now we look at the Cameron-Liebler sets arising from regular ovoids in polar spaces of rank 2.

\begin{example}
	Again we use the notation of Example~\ref{ex:CLps}. If $\P$ is a polar space of rank 4 (so $d=2$), then $\mathcal{M}$ is a regular $(m,1)$-ovoid of $\P'$, which is a polar space of rank 2.
	In this case any $m$-ovoid is a regular $(m,1)$-ovoid, see Remark~\ref{rem:ovoida1}. The $m$-ovoids of classical generalized quadrangles are well studied objects, though there are still some major open problems related to their classification.
	For example, while it is known~\cite{BLP2009} that there are no $m$-ovoids of $\mathcal{H}(4, q^{2})$ when $m \leq \sqrt{q}$, it is an open problem whether any exist for $m > \sqrt{q}$.
	
	Much more is known about $m$-ovoids of $\mathcal{Q}^{-}(5,q)$, which correspond to $m$-covers of $\mathcal{H}(3,q^{2})$.
	Segre~\cite{SegreHemi65} showed that any non-trivial $m$-cover of $\mathcal{H}(3,q^{2})$ for odd $q$ must have $m = \frac{q+1}{2}$, and termed these objects \emph{hemisystems}. He also gave an example of a hemisystem of $\mathcal{H}(3,9)$.
	While this was the only known hemisystem for many years, Cossidente and Penttila~\cite{CossPentt} eventually found an infinite family of hemisystems in $\mathcal{H}(3,q^{2})$ giving examples for all odd $q$,
	and since then a few other infinite families of hemisystems have been found,
	see~\cite{BambergetalHemi, CossPaveseHemi, KorchmarosHemi, SmaldoreHemi}.
	On the other hand, Bruen and Hirschfeld~\cite{BruenHirsch} showed that there are no non-trivial $m$-covers of $\mathcal{H}(3,q^{2})$ (and so no $m$-ovoids of $\mathcal{Q}^{-}(5,q)$) when $q$ is even.  Note that a $\left(\frac{q+1}{2}\right)$-ovoid of $\mathcal{Q}^{-}(5,q)$ does not contain an embedded $\mathcal{Q}(4,q)$, and hence gives rise to a non-trivial Cameron-Liebler set with parameter $\frac{q^{2}(q+1)(q-1)}{2}$ in $\mathcal{Q}(8,q)$.

	There are many examples, and open parameters, for $m$-ovoids of $\mathcal{Q}(4,q)$.
	There are always ovoids ($1$-ovoids) of $\mathcal{Q}(4,q)$ for all $q$ arising from an embedded elliptic quadric; these give rise to trivial Cameron-Liebler sets in $\mathcal{Q}^{+}(7,q)$, see Remark~\ref{rem:newconstructiontrivial}.
	There are also 1-ovoids that are not elliptic quadrics for all odd non-prime $q$, and for all $q = 2^{2h+1} \geq 8$, see~\cite{PWill} for an overview of the known examples.
	Each of these ovoids gives rise to a non-trivial Cameron-Liebler set with parameter $q(q-1)$ in $\mathcal{Q}^{+}(7,q)$.
	Drudge~\cite{DrudgeQ4q} gave a construction of $2$-ovoids, and later Cossidente et al~\cite{CossEtal} showed the existence of $m$-ovoids for all values of $m$, in $\cQ(4,q)$ for $q$ even.
	The latter consist of a union of $2$-ovoids along with an elliptic quadric (if $m$ is odd).
	For odd values of $q$, some sporadic examples found through a computer search were described by Bamberg, Law and Penttila~\cite{BLP2009} for $(q,m) \in \{ (5,2), (7,3), (9,3), (9,4), (11,5) \}$.
	There are also infinite families of $\left(\frac{q-1}{2}\right)$-ovoids known when $q$ is odd (one family for $q \equiv 3 \bmod{4}$ generalizes some of the above mentioned sporadic examples~\cite{FMXOv}, while the other known family~\cite{FengTao} for $q \equiv 1\bmod{4}$ does not).
	Finally, if we have an $m$-ovoid of $\mathcal{Q}^{-}(5,q)$ and take an intersection with a $\mathcal{Q}(4,q)$, then we get an $m$-ovoid of $\mathcal{Q}(4,q)$; thus we have $\frac{q+1}{2}$-ovoids of $\mathcal{Q}(4,q)$ for all odd $q$.
	Each of these $m$-ovoids with $m\geq2$ gives rise to a non-trivial example of a Cameron-Liebler set of generators in $\cQ^+(7,q)$.
\end{example}

Finally we give an example of a Cameron-Liebler set arising through a regular ovoid in a polar space of rank 3.

\begin{example}
	Consider a parabolic quadric $\mathcal{Q}\cong \mathcal{Q}(6,q)$, $q$ odd,
	with an embedded $\mathcal{Q}^{-}\cong\mathcal{Q}^{-}(5,q)$. Let $H$ be a hemisystem, that is a $\left(\frac{q+1}{2}\right)$-ovoid, of $\mathcal{Q}^{-}$.
	Now, let $M$ be the set of all lines of $\mathcal{Q}$ that  meet $\mathcal{Q}^{-}$ in a point of $H$ but are not contained in it.
	Then $M$ is a regular $\left(\frac{q(q+1)}{2},2\right)$-ovoid. We can check this directly using Definition~\ref{def:regularovoida}.
	It is clear that the number of lines of $M$ in relation $R^{1,2}_{2,2}$ with a fixed line of $M$ is constant.
	If $\ell$ is a line not in $M$, then it either is a line of $\mathcal{Q}^{-}$ or else a line not in $\mathcal{Q}^{-}$, but meeting $\mathcal{Q}^{-}$ in a point not in $H$.
	In the former case the number of lines of $M$ in relation $R^{1,2}_{2,2}$ with $\ell$ equals
	\[
	\frac{q+1}{2}\left((q+1)\left(q^{2}+1\right)-q(q+1)-1-q^{2}\right)=\frac{1}{2}q^{2}(q+1)(q-1)
	\]
	since through each point of $H$ on $\ell$ all lines must be counted, except the ones in a plane through $\ell$ and the ones in $\mathcal{Q}^{-}$.
	In the latter case the number of lines of $M$ in relation $R^{1,2}_{2,2}$ with $\ell$ equals
	\[
	q\left(\frac{1}{2}(q+1)\left(q^{2}+1\right)-(q+1)\frac{q+1}{2}\right)=\frac{1}{2}q^{2}(q+1)(q-1)
	\]
	since through each of the points of $\ell$ outside $\mathcal{Q}^{-}$ we must count all lines through a point of $H$, except the ones in a plane through $\ell$.
	As the number is the same in both cases, $M$ is indeed a regular ovoid.
	\par Using this regular ovoid in $\mathcal{Q}(6,q)$ in Example~\ref{ex:CLps}, $q$ odd, we obtain a non-trivial Cameron-Liebler set in $\mathcal{Q}^{+}(9,q)$, $q$ odd, with parameter $\frac{q^{2}(q+1)(q-1)}{2}$.
\end{example}

	
	\subsection*{Acknowledgements}
	Maarten De Boeck and Jozefien D'haeseleer were partially supported by the Croatian Science Foundation under the project 5713.\\
	Jozefien D’haeseleer is supported by the Research Foundation Flanders (FWO) through the grant 1218522N.\\
	Morgan Rodgers is supported by the SFB-TRR 195 `Symbolic Tools in Mathematics and their Application' of the German Research Foundation (DFG).
	
	\bibliographystyle{abbrv}
	\bibliography{main.bib}

	\newpage
	\appendix
	\renewcommand{\thesection}{\Alph{section}}

\section{The eigenvalue calculations}\label{ap:calculations}

Recall that we introduced $\psi_{d,r,i,s,k}$ and $\chi_{d,i,r,n,t,l}$ before, in Definition~\ref{def:psichi}. The first aim is to give a general formula for $\psi_{d,r,i,s,k}$. We will use the following results.

\begin{lemma}[{\cite[p.65 and Lemma  4.3.4(ii)]{phdvanhove}}]\label{lem:psirecursive}
	For integers $d,r,i$ with $0\leq r\leq d-i$ and $0\leq i\leq r$ we have $\psi_{d,r,i,r,d-i}=1$ and for integers $d,r,i,s,k$ with $1\leq s+1\leq r\leq d-i\leq k\leq\min(d,d-i+r-s)-1$ and $0\leq i\leq r$, we have
	\begin{multline*}
		-\gs{d - i + r - s - k}{1}{q}\psi_{d,r,i,s,k}  =q^{r-s-1}\gs{s+1}{1}{q}\psi_{d,r,i,s+1,k} \\
		+q^{d-i+r-s-k-1}\gs{k -d + i +1}{1}{q}\psi_{d,r,i,s,k+1}\:.
	\end{multline*}
\end{lemma}

Explicit formulas for $\psi_{d,r,i,s,k}$ were proved for the cases $s+k=d-i+r$, $k=d$ and $s=0$ in~\cite[Lemma 4.3.5]{phdvanhove} and for the cases $i \in \{ 0,1 \}$ in~\cite[Lemma 11]{metsch}. We now prove the following.

\begin{lemma}
	For integers $d,r,i,j,t$ with $0\leq r\leq d-i$, $0\leq j\leq i\leq r$ and $0\leq t\leq r-j$ we have that
	\begin{multline*}
		\psi_{d,r,i,r-j-t,d-i+j} =
		{(-1)}^{j+t}q^{\binom{j}{2}+\binom{t}{2}+it}\gs{i}{j}{q}\prod_{m=1}^{j}(q^{i-m+e}+1) 
		\cdot \sum_{\ell=0}^{t}{(-1)}^{\ell}\gs{i-j}{\ell}{q}\gs{r-i}{t-\ell}{q}q^{\ell(e-1-t+\ell)}\:.
	\end{multline*}
\end{lemma}
\begin{proof}
	We prove this result using induction on $t$. The case $t=0$ follows from~\cite[Lemma 4.3.5(i)]{phdvanhove}. For the induction step, going from $t$ to $t+1$, we apply the formula from Lemma~\ref{lem:psirecursive} for $(s,k)=(r-j-t-1,d-i+j)$. We find
	{\scriptsize\begin{multline*}
			-\gs{t+1}{1}{q}\psi_{d,r,i,r-j-(t+1),d-i+j} \\
			\begin{aligned}
				& = q^{j+t}\gs{r-j-t}{1}{q}\psi_{d,r,i,r-j-t,d-i+j}+q^{t}\gs{j+1}{1}{q}\psi_{d,r,i,r-t-(j+1),d-i+(j+1)}                                                                                                                \\
				& = \begin{multlined}[t]
					{(-1)}^{j+t}q^{\binom{j+1}{2} +\binom{t+1}{2}+it} \\
					\cdot \prod_{m=1}^{j}(q^{i-m+e}+1)
					\left(\sum_{\ell=0}^{t}{(-1)}^{\ell}\gs{r-j-t}{1}{q}\gs{i}{j}{q}\gs{i-j}{\ell}{q}\gs{r-i}{t-\ell}{q}q^{\ell(e-1-t+\ell)} \right. \\
					\left. -\sum_{\ell=0}^{t}{(-1)}^{\ell}(q^{i-j-1+e}+1)
					\cdot \gs{j+1}{1}{q}\gs{i}{j+1}{q}\gs{i-j-1}{\ell}{q}\gs{r-i}{t-\ell}{q}q^{\ell(e-1-t+\ell)}\right)
				\end{multlined}                                                             \\
				& = \begin{multlined}[t]
					{(-1)}^{j+t}q^{\binom{j+1}{2}+\binom{t+1}{2}+it} \\
					\cdot \gs{i}{j}{q}\prod_{m=1}^{j}(q^{i-m+e}+1)
					\left(\sum_{\ell=0}^{t}{(-1)}^{\ell}\gs{r-j-t}{1}{q}\gs{i-j}{\ell}{q}\gs{r-i}{t-\ell}{q}q^{\ell(e-1-t+\ell)}\right. \\
					\left.-\sum_{\ell=0}^{t}{(-1)}^{\ell}(q^{i-j-1+e}+1)\gs{i-j}{1}{q}\gs{i-j-1}{\ell}{q}\gs{r-i}{t-\ell}{q}q^{\ell(e-1-t+\ell)}\right)
				\end{multlined}                                                                                    \\
				& = \begin{multlined}[t]
					{(-1)}^{j+t}q^{\binom{j+1}{2}+\binom{t+1}{2}+it} \\
					\cdot \gs{i}{j}{q}\prod_{m=1}^{j}(q^{i-m+e}+1)
					\left(\sum_{\ell=0}^{t}{(-1)}^{\ell}\gs{i-j}{\ell}{q}\gs{r-i}{t-\ell}{q}q^{\ell(e-1-t+\ell)}\left(\gs{r-j-t}{1}{q}-\gs{i-j-\ell}{1}{q}\right)\right.                                                      \\
					\left.-\sum_{\ell=0}^{t}{(-1)}^{\ell}\gs{\ell+1}{1}{q}\gs{i-j}{\ell+1}{q}\gs{r-i}{t-\ell}{q}q^{(\ell+1)(e-1)-\ell(t-\ell)+i-j}\right)
				\end{multlined}        \\
				& = \begin{multlined}[t]
					{(-1)}^{j+t}q^{\binom{j}{2}+\binom{t+1}{2}+i(t+1)}\gs{i}{j}{q}\prod_{m=1}^{j}(q^{i-m+e}+1)                                                                                                                         \\
					\left(\sum_{\ell=0}^{t}{(-1)}^{\ell}\gs{i-j}{\ell}{q}\gs{r-i}{t-\ell}{q}q^{\ell(e-1-t+\ell)-\ell}\gs{r-i-t+\ell}{1}{q}\right.                                                                                 \\
					\left.-\sum_{\ell=1}^{t+1}{(-1)}^{\ell-1}\gs{\ell}{1}{q}\gs{i-j}{\ell}{q}\gs{r-i}{t-\ell+1}{q}q^{\ell(e-1)-(\ell-1)(t-\ell+1)}\right)
				\end{multlined}      
			\end{aligned}
	\end{multline*}}

	{\scriptsize\begin{multline*}	
			\begin{aligned}
				& = \begin{multlined}[t]
					{(-1)}^{j+t}q^{\binom{j}{2}+\binom{t+1}{2}+i(t+1)} \\
					\cdot \gs{i}{j}{q}\prod_{m=1}^{j}(q^{i-m+e}+1)
					\sum_{\ell=0}^{t+1}{(-1)}^{\ell}\gs{i-j}{\ell}{q}\gs{r-i}{t-\ell+1}{q}q^{\ell(e-1-(t+1)+\ell)} \\
					\cdot \left(\gs{t-\ell+1}{1}{q}+\gs{\ell}{1}{q}q^{t-\ell+1}\right)
				\end{multlined}                                                                                                                     \\
				& = \begin{multlined}[t]
					{(-1)}^{j+t}q^{\binom{j}{2}+\binom{t+1}{2}+i(t+1)} \\
					\cdot \gs{i}{j}{q}\prod_{m=1}^{j}(q^{i-m+e}+1)
					\sum_{\ell=0}^{t+1}{(-1)}^{\ell}\gs{i-j}{\ell}{q}\gs{r-i}{t-\ell+1}{q}q^{\ell(e-1-(t+1)+\ell)}\gs{t+1}{1}{q}\:.
				\end{multlined}
			\end{aligned}
	\end{multline*}}

	Dividing both side by $-\gs{t+1}{1}{q}$ the formula for the case $t+1$ follows.
\end{proof}

\begin{corollary}\label{cor:psivalues}
	For integers $d,r,i,s,k$ with $0\leq s\leq r\leq d-i\leq k\leq\min(d,d-i+r-s)$ and $0\leq i\leq r$ we have that
	{\scriptsize\begin{multline*}
			\psi_{d,r,i,s,k} = {(-1)}^{r-s}q^{\binom{k-d+i}{2}+\binom{r-k+d-i-s}{2}+i(r-k+d-i-s)}\gs{i}{d-k}{q} \\
			\cdot \prod_{m=1}^{\crampedclap{k-d+i}}(q^{i-m+e}+1)
			\sum_{\ell=0}^{\crampedclap{r-k+d-i-s}}{(-1)}^{\ell}\gs{d-k}{\ell}{q}\gs{r-i}{k-d+s+\ell}{q}q^{\ell(e-1-r+k-d+i+s+\ell)}\:.
	\end{multline*}}
\end{corollary}

Now, we focus on the values $\chi_{d,i,r,n,t,\ell}$. First we introduce some notation and then we mention a recursive relation for the values $\chi_{d,i,r,n,t,\ell}$. This notation and definition is due to~\cite[Definition 3.7]{eisfeld}.

\begin{definition}\label{def:gamma}
	For integers $a,b,c,s,k,t,\ell$ we set
	\[
	\gamma_{a,b,s,k,c,t,\ell}=\gs{s}{t}{q}\gs{k-a}{\ell-a}{q}\gs{a+b-k-s}{a+c-\ell-t}{q}\:q^{(s-t)(\ell-a)+(k+s-\ell-t)(a+c-\ell-t)}\:.
	\]
\end{definition}

\begin{lemma}\label{lem:chirecursive}
	For integers $d,i,r,t,\ell$ with $0\leq t\leq r\leq \ell \leq \min(d,2r-t)$ and $0\leq i\leq\min(r,d-r)$, we have
	\begin{align*}
		\chi_{d,i,r,r,t,\ell}=\sum_{s=t}^{d-i+t-r}\sum_{k=\max(d-i,\ell)}^{\min(d,d-i+\ell-r+t-s)}\gamma_{r,d-i,s,k,r,t,\ell}\:\psi_{d,r,i,s,k}
	\end{align*}
	For integers $d,i,r,n,t,\ell$ with $0\leq t\leq r\leq n\leq \ell \leq \min(d,r+n-t)$ and $0\leq i\leq\min(r,d-n)$ we have
	\begin{multline*}
		\gs{\ell-r}{1}{q}q^{r+n-t-\ell}\chi_{d,i,r,n,t,\ell}
		= \gs{n-r}{1}{q}\chi_{d,i,r,n-1,t,\ell-1} 
		\shoveright{-q^{n-t-1}\gs{t+1}{1}{q}\chi_{d,i,r,n,t+1,\ell-1}} \quad\qquad{}\\
		-\gs{n+r-t-\ell+1}{1}{q}\chi_{d,i,r,n,t,\ell-1}\:.
	\end{multline*}
\end{lemma}

Recall that $\chi_{d,i,r,r,t,\ell}$ is defined to be zero if the condition $0\leq t\leq r\leq n\leq \ell \leq \min(d,r+n-t)$ or the condition $0\leq i\leq\min(r,d-n)$ is not met.
\par Explicit formulas for $\chi_{d,i,r,r,t,\ell}$ were proved in~\cite[Lemma 4.3.9]{phdvanhove} for the cases $(t,\ell) \in \{ (r-1,n), (r-1,n+1) \}$ with $r\geq1$ (one of them is mentioned in Lemma~\ref{lem:psichivalues}),
and in~\cite[Lemma 12]{metsch} for $(i,\ell) \in \{ (0,n),(1,n) \}$ and for $(i,t,\ell) \in \{ (0,0,n+1),(1,0,n+1) \}$. We mention one of these results.

\begin{lemma}\label{lem:chi_n=l}
	For integers $d,r,n,t$ with $0\leq t\leq r\leq n \leq d-1$ and $1\leq r$ we have that
	\begin{align*}
		\chi_{d,1,r,n,t,n} & ={(-1)}^{r-t}q^{\binom{r-t}{2}+(r-t)(d-n-1)}\left(\gs{r}{t}{q}-(q^{e}+1)\gs{r-1}{t}{q}\right)
	\end{align*}
\end{lemma}

Using the formulas from Corollary~\ref{cor:psivalues} a proof different from the one in~\cite{metsch} can be given.
\par We will now prove an explicit formula for the case $(i,\ell)=(1,n+1)$.

\begin{lemma}\label{lem:chi_n+1=l}
	For integers $d,r,n,t$ with $0\leq t\leq r\leq n \leq \min(d,r+n-t)-1$ and $r\geq 1$ we have that
	{\scriptsize\begin{multline*}
			\chi_{d,1,r,n,t,n+1}  = {(-1)}^{r-t}q^{\binom{r-t}{2}+(r-t-1)(d-n-2)}\gs{r-1}{t}{q} \\
			\cdot \left((q^{e}+1)\left(q^{d-t-2-n+r}+(q^{r-t-1}-1)\gs{d-n-1}{1}{q}\right) \right. 
			\left. -(q^{r}-1)\gs{d-n-1}{1}{q}\right)\:.
	\end{multline*}}
\end{lemma}
\begin{proof}
	We prove this result using induction on $n$. For $n=r$ we derive from the first formula in Lemma~\ref{lem:chirecursive} that
	{\scriptsize\begin{align*}
			\chi_{d,1,r,r,t,r+1} & =\sum_{s=t}^{d-1+t-r}\sum_{k=\max(d-1,r+1)}^{\min(d,d+t-s)}\gamma_{r,d-1,s,k,r,t,r+1}\:\psi_{d,r,1,s,k} \\
			& =
			\begin{cases}
				\gamma_{d-1,d-1,t,d,d-1,t,d}\:\psi_{d,d-1,1,t,d} & r=d-1                                      \\
				\gamma_{r,d-1,t,d-1,r,t,r+1}\:\psi_{d,r,1,t,d-1}+\gamma_{r,d-1,t,d,r,t,r+1}\:\psi_{d,r,1,t,d} \\\quad+\gamma_{r,d-1,t+1,d-1,r,t,r+1}\:\psi_{d,r,1,t+1,d-1}&1\leq r\leq d-2
			\end{cases}
	\end{align*}}
	So, if $1\leq r\leq d-2$ we find using Corollary~\ref{cor:psivalues} and Definition~\ref{def:gamma} that
	{\scriptsize\begin{align*}
			\chi_{d,1,r,r,t,r+1} & =
			\begin{multlined}[t]
				\gs{d-r-1}{1}{q}\gs{r-t}{1}{q}q^{(d-r-2)(r-t-1)}{(-1)}^{r-t}q^{\binom{r-t+1}{2}}\left(\gs{r-1}{t-1}{q}-\gs{r-1}{t}{q}q^{e+t-r}\right)          \\
				+\gs{d-r}{1}{q}q^{(d-r-1)(r-t-1)}{(-1)}^{r-t}q^{\binom{r-t}{2}}\left(q^{e}+1\right)\gs{r-1}{t}{q}                                \\
				+\gs{t+1}{1}{q}\gs{d-r-1}{1}{q}q^{1+(d-r-1)(r-t-1)}{(-1)}^{r-t+1}q^{\binom{r-t}{2}}                                              \\
				\left(\gs{r-1}{t}{q}-\gs{r-1}{t+1}{q}q^{e+t-r+1}\right)
			\end{multlined}                                                  \\
			& =
			\begin{multlined}[t]
				{(-1)}^{r-t}q^{\binom{r-t}{2}+(d-r-2)(r-t-1)}\left(q^{r-t}\gs{d-r-1}{1}{q}\gs{r-t}{1}{q}\gs{r-1}{t-1}{q}\right.                       \\
				-q^{e}\gs{d-r-1}{1}{q}\gs{r-t}{1}{q}\gs{r-1}{t}{q}+\left(q^{e}+1\right)\gs{d-r}{1}{q}q^{r-t-1}\gs{r-1}{t}{q}                     \\
				\left.-q^{r-t}\gs{t+1}{1}{q}\gs{d-r-1}{1}{q}\gs{r-1}{t}{q}+q^{e+1}\gs{t+1}{1}{q}\gs{d-r-1}{1}{q}\gs{r-1}{t+1}{q}\right)
			\end{multlined}                  \\
			& =
			\begin{multlined}[t]
				{(-1)}^{r-t}q^{\binom{r-t}{2}+(d-r-2)(r-t-1)}\left(q^{r-t}\gs{d-r-1}{1}{q}\gs{t}{1}{q}\gs{r-1}{t}{q}\right. \\
				-q^{e}\gs{d-r-1}{1}{q}\gs{r-t}{1}{q}\gs{r-1}{t}{q}+\left(q^{e}+1\right)\gs{d-r}{1}{q}q^{r-t-1}\gs{r-1}{t}{q}                                      \\
				\left.-q^{r-t}\gs{t+1}{1}{q}\gs{d-r-1}{1}{q}\gs{r-1}{t}{q}+q^{e+1}\gs{r-t-1}{1}{q}\gs{d-r-1}{1}{q}\gs{r-1}{t}{q}\right)
			\end{multlined} \\
			& =
			\begin{multlined}[t]
				{(-1)}^{r-t}q^{\binom{r-t}{2}+(d-r-2)(r-t-1)}\gs{r-1}{t}{q}\left(q^{r-t}\gs{d-r-1}{1}{q}\gs{t}{1}{q}\right. \\
				-q^{r-t}\gs{t+1}{1}{q}\gs{d-r-1}{1}{q}+\left(q^{e}+1\right)\gs{d-r}{1}{q}q^{r-t-1}\\
				\left.-q^{e}\gs{d-r-1}{1}{q}\gs{r-t}{1}{q}+q^{e+1}\gs{r-t-1}{1}{q}\gs{d-r-1}{1}{q}\right)
			\end{multlined}                                                                   \\
			& =
			\begin{multlined}[t]
				{(-1)}^{r-t}q^{\binom{r-t}{2}+(d-r-2)(r-t-1)}\gs{r-1}{t}{q}\\
				\left(-q^{r}\gs{d-r-1}{1}{q}+\left(q^{e}+1\right)\gs{d-r}{1}{q}q^{r-t-1}-q^{e}\gs{d-r-1}{1}{q}\right)
			\end{multlined}                                                                  \\
			& =
			\begin{multlined}[t]
				{(-1)}^{r-t}q^{\binom{r-t}{2}+(d-r-2)(r-t-1)}\gs{r-1}{t}{q}\\
				\left(\left(q^{e}+1\right)\left(q^{d-t-2}+\gs{d-r-1}{1}{q}q^{r-t-1}\right)-q^{r}\gs{d-r-1}{1}{q}-q^{e}\gs{d-r-1}{1}{q}\right)\:,
			\end{multlined}
	\end{align*}}
	which simplifies to the formula in the statement of the lemma. If $r=d-1$ we find using Corollary~\ref{cor:psivalues} and Definition~\ref{def:gamma} that
	\begin{align*}
		\chi_{d,1,d-1,d-1,t,d} & ={(-1)}^{d-1-t}q^{\binom{d-1-t}{2}}\left(q^{e}+1\right)\gs{d-2}{t}{q}\:,
	\end{align*}
	which matches the formula in the statement of the lemma. This concludes the case $n=r$.

	For the induction step, going from $n$ to $n+1$, we use the recursive relation from Lemma~\ref{lem:chirecursive} and the result from Lemma~\ref{lem:chi_n=l}.
	{\scriptsize\begin{align*}
			\chi_{d,1,r,n+1,t,n+2} & =
			\begin{multlined}[t]
				\frac{q^{t+1-r}}{\gs{n+2-r}{1}{q}}\left(\gs{n+1-r}{1}{q}\chi_{d,1,r,n,t,n+1} \right. \\
				\left. -\gs{t+1}{1}{q}q^{n-t}\chi_{d,1,r,n+1,t+1,n+1}
				-\gs{r-t}{1}{q}\chi_{d,1,r,n+1,t,n+1}\right)
			\end{multlined}                                                                                                            \\
			& =
			\begin{multlined}[t]
				\frac{q^{t+1-r}}{\gs{n+2-r}{1}{q}}\left(\gs{n+1-r}{1}{q}{(-1)}^{r-t}q^{\binom{r-t}{2}+(r-t-1)(d-n-2)}\gs{r-1}{t}{q}\right.\\
				\left.\left[(q^{e}+1)\left(q^{d-t-2-n+r}+(q^{r-t-1}-1)\gs{d-n-1}{1}{q}\right)-(q^{r}-1)\gs{d-n-1}{1}{q}\right]\right.\\
				-\gs{t+1}{1}{q}q^{n-t}{(-1)}^{r-t-1}q^{\binom{r-t-1}{2}+(r-t-1)(d-n-2)}\left(\gs{r}{t+1}{q}-(q^{e}+1)\gs{r-1}{t+1}{q}\right)\\
				\left.-\gs{r-t}{1}{q}{(-1)}^{r-t}q^{\binom{r-t}{2}+(r-t)(d-n-2)}\left(\gs{r}{t}{q}-(q^{e}+1)\gs{r-1}{t}{q}\right)\right)
			\end{multlined}  \\
			& =
			\begin{multlined}[t]
				\frac{q^{t+1-r}}{\gs{n+2-r}{1}{q}}{(-1)}^{r-t}q^{\binom{r-t}{2}+(r-t-1)(d-n-2)}\left(\gs{n+1-r}{1}{q}\gs{r-1}{t}{q}\right.\\
				\left.\left[(q^{e}+1)\left(q^{d-t-2-n+r}+(q^{r-t-1}-1)\gs{d-n-1}{1}{q}\right)-(q^{r}-1)\gs{d-n-1}{1}{q}\right]\right.\\
				+q^{n-r+1}\left(\gs{r-1}{t}{q}\gs{r}{1}{q}-(q^{e}+1)\gs{r-1}{t}{q}\gs{r-t-1}{1}{q}\right)\\
				\left.-q^{d-n-2}\left(\gs{r-1}{t}{q}\gs{r}{1}{q}-(q^{e}+1)\gs{r-1}{t}{q}\gs{r-t}{1}{q}\right)\right)
			\end{multlined} \\
			& =
			\begin{multlined}[t]
				\frac{1}{\gs{n+2-r}{1}{q}}{(-1)}^{r-t}q^{\binom{r-t}{2}+(r-t-1)(d-n-3)}\gs{r-1}{t}{q}\left((q^{e}+1)\vphantom{\left(\gs{n+1-r}{1}{q}\right)}\right. \\
				\left(q^{d-t-2-n+r}\gs{n+1-r}{1}{q}+(q^{r-t-1}-1)\gs{d-n-1}{1}{q}\gs{n+1-r}{1}{q}\right.\\
				\left.\qquad\qquad-q^{n-r+1}\gs{r-t-1}{1}{q}+q^{d-n-2}\gs{r-t}{1}{q}\right)\\
				-(q^{r}-1)\gs{d-n-1}{1}{q}\gs{n+1-r}{1}{q}+q^{n-r+1}\gs{r}{1}{q}\left.-q^{d-n-2}\gs{r}{1}{q}\right)
			\end{multlined}                                                                       \\
			& =
			\begin{multlined}[t]
				\frac{1}{\gs{n+2-r}{1}{q}}{(-1)}^{r-t}q^{\binom{r-t}{2}+(r-t-1)(d-n-3)} \gs{r-1}{t}{q} \left((q^{e}+1) \vphantom{\left(\gs{n+1-r}{1}{q}\right)}\right. \\
				\left(q^{d-t-3-n+r}\gs{n+2-r}{1}{q}-q^{n-r+1}\gs{r-t-1}{1}{q}+q^{d-n-2}\gs{r-t-1}{1}{q}\right.\\
				\left.+(q^{r-t-1}-1) \left(\gs{d-n-2}{1}{q}\gs{n+2-r}{1}{q} +\frac{q^{n+1-r}-q^{d-n-2}}{q-1}\right)\right)\\
				\left.-(q^{r}-1)\left(\gs{d-n-1}{1}{q}\gs{n+1-r}{1}{q} -\frac{q^{n+1-r}-q^{d-n-2}}{q-1}\right)\right)
			\end{multlined}                                                            \\
			& =
			\begin{multlined}[t]
				{(-1)}^{r-t}q^{\binom{r-t}{2}+(r-t-1)(d-n-3)}\gs{r-1}{t}{q}\\
				\left((q^{e}+1)\left(q^{d-t-3-n+r}+(q^{r-t-1}-1)\gs{d-n-2}{1}{q}\right)-(q^{r}-1)\gs{d-n-2}{1}{q}\right)
			\end{multlined}
	\end{align*}}
	This completes the proof of the induction step.
\end{proof}

The importance of the values $\chi_{d,i,r,n,t,\ell}$, which were introduced as a scaling factor between eigenvalues (see Definition~\ref{def:psichi}), is that they are directly connected to the general eigenvalues.

\begin{theorem}[{\cite[Theorem 3.8]{eisfeld}}]\label{th:mu}
	The eigenvalue of $C^{s,k}_{n,n}$ on $V^{n}_{r,i}$ is
	\[
	\mu_{d,n,s,k,r,i}= \sum^{s}_{t=0}\sum^{\min\{k,n+r-t\}}_{\ell=s+r-t} \alpha_{d,(n,r),(t,\ell),n,(s,k)}\:\chi_{d,i,r,n,t,\ell}\:.
	\]
\end{theorem}

We use this theorem to compute the eigenvalues $\mu_{d,d-1,s,k,r,1}$ of the relations $R^{s,k}_{d-1,d-1}$ on the subspaces $V^{d-1}_{r,1}$.

\begin{theorem}\label{th:cor:eigenval}
	For integers $d,r,s,k$ with $1\leq r\leq d-1$, $0\leq s\leq d-1$ and $d-1\leq k\leq\min(d,2d-s-2)$ we have that the eigenvalue of $C^{s,k}_{d-1,d-1}$ on $V^{d-1}_{r,1}$ equals
	\[
	\mu_{d,d-1,s,k,r,1}=q^{e(d-r-s-1)}\sum^{s}_{t=0} \gs{d-r-1}{s-t}{q}{(-1)}^{r-t}q^{\binom{r-t}{2}+\binom{d-r-s+t+1}{2}+et-1}\:\mu'_{d,k,r,t}
	\]
	with
	$\mu'_{d,k,r,t}=
	\begin{cases}
		\left(\gs{r}{t}{q}-(q^{e}+1)\gs{r-1}{t}{q}\right) & k=d-1 \\
		(q^{e}+1)\gs{r-1}{t}{q}                           & k=d
	\end{cases}$
\end{theorem}
\begin{proof}
	By Theorem~\ref{th:mu} and since $\chi_{d,1,r,d-1,t,\ell}=0$ if $\ell<d-1$ we have that
	\begin{equation*}
		\mu_{d,d-1,s,k,r,1}= \sum^{s}_{t=0}
		\quad\qquad\sum^{\mathclap{\min\{k,d+r-t-1\}}}_{\mathclap{\ell=\max\{s+r-t,d-1\}}} \alpha_{d,(d-1,r),(t,\ell),d-1,(s,k)}\cdot\chi_{d,1,r,d-1,t,\ell}\:.
	\end{equation*}
	For $k=d-1$ we find, using Lemmas~\ref{lem:alphavalues} and~\ref{lem:chi_n=l}, that
	{\scriptsize\begin{align*}
			\mu_{d,d-1,s,d-1,r,1}
			& =
			\sum^{\min\{s,r\}}_{\mathclap{t=\max\{0,s+r-d+1\}}} \alpha_{d,(d-1,r),(t,d-1),d-1,(s,d-1)}\cdot\chi_{d,1,r,d-1,t,d-1} \\
			& =
			\begin{multlined}[t]
				\sum^{\min\{s,r\}}_{\mathclap{t=\max\{0,s+r-d+1\}}} \gs{d-r-1}{s-t}{q}q^{\frac{1}{2}(d-1-r+t-s)(d-r+t-s+2e+2)}\\
				\cdot{(-1)}^{r-t}q^{\binom{r-t}{2}} \left(\gs{r}{t}{q}-(q^{e}+1)\gs{r-1}{t}{q}\right)
			\end{multlined}         \\
			& =
			\begin{multlined}[t]
				\sum^{s}_{t=0} \gs{d-r-1}{s-t}{q}{(-1)}^{r-t}q^{\binom{r-t}{2}+\binom{d-r-s+t}{2}+(e+1)(d-r-s+t-1)}\\
				\left(\gs{r}{t}{q}-(q^{e}+1)\gs{r-1}{t}{q}\right)
			\end{multlined}                    \\
			& =
			\begin{multlined}[t]
				q^{e(d-r-s-1)}\sum^{s}_{t=0} \gs{d-r-1}{s-t}{q}{(-1)}^{r-t}q^{\binom{r-t}{2}+\binom{d-r-s+t+1}{2}+et-1}\\
				\left(\gs{r}{t}{q}-(q^{e}+1)\gs{r-1}{t}{q}\right)
\:.
			\end{multlined}
	\end{align*}}
	
	For $k=d$ we find, using Lemmas~\ref{lem:alphavalues},~\ref{lem:chi_n=l} and~\ref{lem:chi_n+1=l}, that
	{\scriptsize\begin{align*}
			\mu_{d,d-1,s,d,r,1}
			& =
			\begin{multlined}[t]
				\sum^{\min\{s,r\}}_{\mathclap{t=\max\{0,s+r-d+1\}}} \alpha_{d,(d-1,r),(t,d-1),d-1,(s,d)} \cdot\chi_{d,1,r,d-1,t,d-1}\\
				+\sum^{\min\{s,r-1\}}_{\mathclap{t=\max\{0,s+r-d\}}} \alpha_{d,(d-1,r),(t,d),d-1,(s,d)} \cdot\chi_{d,1,r,d-1,t,d}
			\end{multlined}                                               \\
			& =
			\begin{multlined}[t]
				\sum^{\min\{s,r\}}_{t=\max\{0,s+r-d+1\}} \gs{d-1-r}{s-t}{q}q^{d-1-r+t-s+\frac{1}{2}(d-2-r+t-s)(d-r+t-s+2e-1)}\\
				\cdot(q^{e}+1)\gs{d-1-r+t-s}{1}{q}{(-1)}^{r-t}q^{\binom{r-t}{2}}\left(\gs{r}{t}{q}-(q^{e}+1)\gs{r-1}{t}{q}\right)\\
				+\sum^{\min\{s,r-1\}}_{t=\max\{0,s+r-d\}}\gs{d-r}{s-t}{q}q^{\frac{1}{2}(d-1-r+t-s)(d-r+t-s+2e)}\gs{d-r+t-s}{1}{q}\\
				\cdot{(-1)}^{r-t}q^{\binom{r-t}{2}}\gs{r-1}{t}{q}(q^{e}+1)
			\end{multlined}                                                  \\
			& =
			\begin{multlined}[t]
				\sum^{s}_{t=0}{(-1)}^{r-t}(q^{e}+1)q^{\binom{r-t}{2}+\binom{d-r-s+t}{2}+e(d-r-s+t-2)}\\
				\cdot\gs{d-1-r}{s-t}{q}\gs{d-1-r+t-s}{1}{q}
				\left(q^{r-t}\gs{r-1}{t-1}{q}-q^{e}\gs{r-1}{t}{q}\right)\\
				+\sum^{s}_{t=0}{(-1)}^{r-t}(q^{e}+1)q^{\binom{r-t}{2}+\binom{d-r-s+t}{2}+e(d-r-s+t-1)}\gs{d-1-r}{s-t}{q}\gs{d-r}{1}{q}\gs{r-1}{t}{q}
			\end{multlined}                               \\
			& =
			\begin{multlined}[t]
				\sum^{s}_{t=0}{(-1)}^{r-t}(q^{e}+1)q^{\binom{r-t}{2}+\binom{d-r-s+t}{2}+e(d-r-s+t-2)}\gs{d-1-r}{s-t}{q}\\
				\left(\gs{d-1-r+t-s}{1}{q}\left(q^{r-t}\gs{r-1}{t-1}{q}-q^{e}\gs{r-1}{t}{q}\right) +\gs{d-r}{1}{q}q^{e} \gs{r-1}{t}{q}\right)
			\end{multlined}                                                  \\
			& =
			\begin{multlined}[t]
				\sum^{s}_{t=0}{(-1)}^{r-t}(q^{e}+1)q^{\binom{r-t}{2}+\binom{d-r-s+t}{2}+e(d-r-s+t-2)}\gs{d-1-r}{s-t}{q}\\
				\left(q^{r-t}\gs{d-1-r+t-s}{1}{q}\gs{r-1}{t-1}{q}+\gs{s-t+1}{1}{q}q^{e+d-1-r+t-s}\gs{r-1}{t}{q}\right)
			\end{multlined}                                                            \\
			& =
			\begin{multlined}[t]
				\sum^{s}_{t=0}{(-1)}^{r-t}(q^{e}+1) q^{\binom{r-t+1}{2} +\binom{d-r-s+t}{2}+e(d-r-s+t-2)} \gs{d-1-r}{s-t+1}{q} \gs{s-t+1}{1}{q} \gs{r-1}{t-1}{q}\\
				+\sum^{s}_{t=0}{(-1)}^{r-t}(q^{e}+1) q^{\binom{r-t}{2} +\binom{d-r-s+t+1}{2}+e(d-r-s+t-1)-1} \gs{d-1-r}{s-t}{q}\\
				\left(\gs{s-t+1}{1}{q}\gs{r-1}{t}{q}\right)
			\end{multlined}       \\
			& =
			\begin{multlined}[t]
				\sum^{s-1}_{t=0}{(-1)}^{r-t-1}(q^{e}+1)q^{\binom{r-t}{2}+\binom{d-r-s+t+1}{2}+e(d-r-s+t-1)}\gs{d-1-r}{s-t}{q}\gs{s-t}{1}{q}\gs{r-1}{t}{q}\\
				+\sum^{s}_{t=0}{(-1)}^{r-t}(q^{e}+1)q^{\binom{r-t}{2}+\binom{d-r-s+t+1}{2}+e(d-r-s+t-1)-1}\gs{d-1-r}{s-t}{q}\\
				\left(\gs{s-t+1}{1}{q}\gs{r-1}{t}{q}\right)
			\end{multlined}                          \\
			& =
			\begin{multlined}[t]
				(q^{e}+1)q^{e(d-r-s-1)} \left(-\sum^{s-1}_{t=0} {(-1)}^{r-t-1} q^{\binom{r-t}{2} +\binom{d-r-s+t+1}{2} +et-1} \gs{d-1-r}{s-t}{q}\gs{r-1}{t}{q}\right. \\
				+{(-1)}^{r-s}q^{\binom{r-s}{2}+\binom{d-r+1}{2}+es-1}\left. \gs{r-1}{s}{q}\right)                       \end{multlined} \\
			& =(q^{e}+1)q^{e(d-r-s-1)}\sum^{s}_{t=0}{(-1)}^{r-t}q^{\binom{r-t}{2}+\binom{d-r-s+t+1}{2}+et-1}\gs{d-1-r}{s-t}{q}\gs{r-1}{t}{q}\:.\qedhere
	\end{align*}}
\end{proof}

\end{document}